\documentclass[a4paper, 10pt]{amsart}
\usepackage{graphicx}
\usepackage{enumerate}
\usepackage{xspace}
\usepackage{boxedminipage}
\usepackage[ruled]{algorithm}
\usepackage{subfigure}
\usepackage{booktabs}
\usepackage{tabularx}
\usepackage[square,sort,compress]{natbib}
\usepackage{hyperref} 
\usepackage{times}
\usepackage{amsmath,amssymb,amsxtra,amsthm}
\usepackage{subfigure}
\usepackage{geometry}                
\usepackage{tikz}
\usepackage{multirow}

\swapnumbers
\theoremstyle{definition}
\numberwithin{equation}{section}

\newtheorem{Pbm}[subsection]{Problem}

\newtheorem{Rem}[subsection]{Remark}

\newcommand{ \bl}{\color{blue}}

\definecolor{MyGreen}{rgb} {0.05,0.4,0.05}

\usepackage{ifthen}
\provideboolean{showchanges}
\setboolean{showchanges}{false}
\newcommand{\changes}[1]{
  \ifthenelse{\boolean{showchanges}} {{\bl{#1}}} {#1}
}
\provideboolean{shownotes}
\setboolean{shownotes}{true}
\newcommand{\standout}[1]{\colorbox{red}{\textcolor{white}{#1}}}
\newcounter{margnote}[page]
\newcommand{\margnotemark}{{\standout{\footnotesize\upshape\texttt{\arabic{margnote}}}}}
\newcommand{\margnote}[2][]{
  \ifthenelse{
    \boolean{shownotes}
  }{\stepcounter{margnote}\margnotemark\marginpar{
      \texttt{
        \begin{minipage}{2cm}
          \raggedright\tiny
          \margnotemark{#1}: 
          #2
        \end{minipage}
  }}}{}
}
\usepackage{upgreek}
\newcommand{\lap}{\ensuremath{\Updelta}}

\renewcommand{\vec}[1]{\ensuremath{\boldsymbol{#1}}}

\newcommand{\pdt}{\ensuremath{\partial_t}}
\newcommand{\mdt}{\ensuremath{\partial^{\bullet}}}

\newcommand{\Reals}{\ensuremath{{\mathbb{R}}}}

\newcommand{\diff}{\ensuremath{{\operatorname{d}}}}
\renewcommand{\O}{\ensuremath{{\Omega}}}

\newcommand{\normal}{\ensuremath{{\vec{\nu}}}}

\newcommand{\G}{\ensuremath{{\Gamma}}}
\newcommand{\Gc}{\ensuremath{{\hat{\G}}}}
\newcommand{\Gt}{\ensuremath{{\G(t)}}}
\newcommand{\Gct}{\ensuremath{{\hat{\G}(t)}}}
\newcommand{\Lp}[1]{\ensuremath{\operatorname{L}_{#1}}}

\newcommand{\ltwon}[2]{\ensuremath{\left\|#1\right\|}_{\Lp{2}#2}}

\newcommand{\lv}{\ensuremath{\left\vert}}
\newcommand{\rv}{\ensuremath{\right\vert}}
\usepackage{upgreek}

\usepackage{color}
\definecolor{MyGreen}{rgb} {0.05,0.4,0.05}
\definecolor{RedViolet}{rgb} {0.1,0.1,0.75}

\usepackage[british]{babel}
\renewenvironment{itemize}{\begin{list}{\labelitemi}{\leftmargin=1em} }
{\end{list}
}

\def\mmm{\mathcal}

\def\bpmat{\begin{pmatrix}}
\def\epmat{\end{pmatrix}}
\def\mathref#1{\ifmmode\mathrm{(\ref{#1})}\else{\rm(\ref{#1})}\fi} 
\def\nref#1{\ifmmode\mathrm{\ref{#1}}\else{\rm\ref{#1}}\fi}

\newcommand \beq{\begin{equation}}
\newcommand \eeq{\end{equation}}
\usepackage{hyperref}
\begin{document}
\bibliographystyle{abbrvnat}
  \title{Parameter identification problems in the modelling of cell motility}
   \date{\today}
\author{Wayne Croft}  
\address[W.~Croft]{School of Life Sciences\\ Queens Medical Centre\\ University of Nottingham\\ Nottingham\\ UK\\ NG7 2UH} 
\author{Charles M Elliott}
\address[C.~M.~Elliott\\B.~Stinner]{Mathematics Institute\\
 Zeeman Building\\ 
 University of Warwick\\ Coventry\\ UK\\
 CV4 7AL}
\author{Graham Ladds}
\address[G.~Ladds\\ C.~Weston]{Division of Biomedical Cell Biology\\ Warwick Medical School\\ University of Warwick\\ Coventry\\ UK\\ CV4 7AL}
\author{Bj\"orn Stinner}
\author{Chandrasekhar Venkataraman}
\address[C.~Venkataraman]{Department of mathematics\\
	    School of Mathematical and Physical Sciences\\
	    University of Sussex\\
	    Falmer\\
	    near Brighton\\
	    UK, BN1 9RF}
\email[C.~Venkataraman]{c.venkataraman@sussex.ac.uk}
\author{Cathryn Weston}
\begin{abstract}
We present a novel parameter identification algorithm for the estimation of parameters in models of cell motility using imaging data of migrating cells. Two alternative formulations of the objective functional that measures the difference between the computed and observed data are proposed and the parameter identification problem is formulated as a minimisation problem of nonlinear least squares type. A Levenberg-Marquardt based optimisation method is applied to the solution of the minimisation problem and the details of the implementation are discussed. A number of numerical experiments are presented which illustrate the robustness of the algorithm to parameter identification in the presence of large deformations and noisy data and parameter identification in three dimensional models of cell motility. An application to experimental data is also presented in which we seek to identify parameters in a model for the monopolar growth of fission yeast cells using experimental imaging data.  
\end{abstract}
\maketitle
 \section{Introduction}\label{sec:intro}
Mathematical modelling and numerical simulation of the directional motility of cells  is of much importance, in part due to the central role directed cell migration plays in several biological phenomena such as embryonic development, cancer, tissue development and immune responses \citep{bray2001cell}. In \citet{venk11chemotaxis} we developed a general framework for the modelling and simulation of cell motility and chemotaxis. The modelling framework we developed is a phenomenological one, in which we propose a  geometric evolution law for the   cell membrane dynamics which accounts for its mechanical properties. For the polarisation of the cell we postulate a reaction diffusion system for species located on the moving cell membrane. Protrusion is then achieved by back-coupling these surface quantities to the geometric equation for the membrane position. We developed an efficient and robust numerical method for the simulation of the model equations and used the method to compare the behaviour of the simulated  cells with experimental data from biological cells migrating in vitro. 

Since the models we proposed are phenomenological there are many parameters in the model equations which do not, as yet, have direct biological counterparts and can thus not be selected from experimental data. Moreover, it remains difficult experimentally to  quantify the forces associated with motility  and only recently has experimental progress been made in this direction \citep{del2007spatio, lombardi2007traction}. Thus even when parameters in the models correspond directly to a physically meaningful quantity, experimental
measurements may be unavailable. Due to rapid advancements in fluorescence microscopy and experimental techniques, the  imaging of migrating cells both in vitro and in vivo is a burgeoning  research field.  High resolution (both temporal and spatial) three dimensional data of migrating cells with concentrations of motility related species tagged with fluorescent marker proteins is available, see \citet{pittet2011intravital} for a review. This information may be used to estimate, otherwise inaccessible, parameters in mathematical models of cell motility and may even be used in formulating the model itself \citep{sbalzarini2013modeling}. This provides the motivation for this study  in which we present an algorithm for the identification of parameters in models of cell motility, or more generally in coupled geometric evolution law - surface partial differential equation models, where (potentially noisy) observations of the position of the cell membrane and associated concentrations of membrane resident species are available at a series of discrete times.  

Parameter identification for semilinear reaction-diffusion systems has received some attention in the literature  \citep{friedman1992parameter, garvie2010efficient,  jiang2000parameter, ackleh1998numerical} and  there have also been some studies focussed on optimal control and related inverse problems in the context of geometric evolution equations and free boundary problems \citep{deckelnick2009optimal,deckelnick2011numerical, hausser2010influence,hausser2012control}.  There has also been some work on parameter estimation in models of cell motility. \citet{satulovsky2008exploring} propose a discrete  rule-based model for cell motility and  optimise the parameters in their model such that their simulated cells have similar motility related statistics, such as persistence length and centroid velocity, to experimentally observed cells.  \citet{milutinovicparameters} consider a simple model for the motion of the cell centroid and use a data assimilation approach to fit the parameters in their model to experimental observations of migrating cells in the embryo of the Zebrafish. To our best knowledge very few works treat parameter estimation in the setting considered in this study   where a geometric evolution law, or more generally a free boundary problem, is coupled to a reaction diffusion system. One notable exception is \citet{hogea2008image} where a parameter identification method based on optimal control theory is proposed for a coupled reaction-diffusion system - free boundary problem modelling the growth of gliomas.

As mentioned above, in \citet{venk11chemotaxis} we have developed an efficient and robust solver for reaction-diffusion system - evolution law based models of cell motility based on a parametric surface finite element method \citep{dziuk2007finite,dziuk2008computational,barrett2007parametric,barrett2008parametric}. The main focus of this study is to present an algorithm for the identification of parameters in the model equations, in which the forward problem may be solved with this parametric surface finite element method. Thus we wish to avoid the consideration of embedded methods such as phase field or level set methods for the solution of the forward problem as these methods may become prohibitive in terms of computational time for large-scale iterative optimisation problems in three dimensional models of cell motility. 

A significant contribution of this work is the construction of appropriate objective functionals for the comparison of the simulated cells with the observations. We propose two versions of the objective functional, both formulated in the continuous setting. The first is based on the sharp interface description of the cell membrane and consists of  the Haussdorff distance between the observed data and the simulated cells together with the difference in concentration between a point on one surface and the concentration at the closest point on the other surface. The second formulation of the objective functional is based on a phase field representation of the surfaces and an extension of the concentrations from the sharp interface to a diffuse interface. This phase field formulation of the objective functional may be computed from the sharp interface description at the times at which the objective is evaluated and thus the parametric finite element method previously developed may still be used for the solution of the model equations. The reasons behind presenting the two different formulations of the objective functional is that the first appears more robust computationally, in certain settings, and is cheaper to compute while the second, unlike the first, is smooth and thus it may be possible to analyse the optimisation algorithm we propose, which assumes smoothness of the objective functional. 

The parameter identification problem is then formulated as a nonlinear least squares problem and we present an iterative optimisation method for its minimisation. The Levenberg-Marquardt method \citep{marquardt1963algorithm}  is considered for the optimisation,  due to the fact that it is a widely used numerical method for nonlinear least squares problems that has proved robust in practice and converges in settings where other standard methods such as Gauss-Newton fail (e.g., when the Jacobian matrix does not have full column rank).   

We discuss  the implementation of the optimisation method and illustrate that the formulation of the objective functional may be adapted to deal with the types of data typically generated in experiments. The performance of the algorithm is examined in various numerical tests,  where we use artificial data generated by numerical simulations. To investigate the robustness of the algorithm, we conduct simulations with noisy data.  We also report on the sensitivity of the two different formulations of the objective functional to changes in the relative contribution of the two distinct sources of error, i.e., the error due to position and the error due to concentration. Simulations where we examine the effect of varying the size of the interfacial width parameter associated with the phase field method are also presented. To illustrate that the algorithm is applicable to the three dimensional models of cell motility, we report on a parameter identification experiment in three dimensions. 

As a proof of concept that the algorithm can be applied to real data, we use the identification algorithm to identify parameters in a  simple model for the monopolar growth of fission yeast cells where the data consists of experimental observations. We describe the experimental setup employed to generate the experimental data and then propose a simple model for the observed monopolar growth in which two parameters are to be identified. A number of numerical simulations that illustrate that the algorithm is able to identify parameters with experimental data are presented and the parameters identified show good agreement with experimental estimates of related parameters. We also show that by changing the relative contribution of the error due to position and concentration (through changing weights) to the  objective functional, it is possible to identify parameters that appear to prioritise the fit to observations of either the positions or the concentrations. This we believe will be of interest to practitioners who have a specific application in mind where the fit to either position or concentration is of more importance.

The remainder of our discussion proceeds as follows. In \S \ref{sec:setup} we define our notation, formulate the forward problem, present the two different formulations of the objective functional we consider and state the identification problem. In \S \ref{sec:optimisation} we discuss optimisation methods for the solution of the identification problem and present  the Levenberg-Marquardt method applied to the solution of the identification problem we consider in this study. In \S \ref{sec:implementation} we discuss the implementation of the algorithm and the application of the algorithm to observations of the form typically generated in experiments. In \S \ref{sec:examples} we report on a series of numerical experiments with artificial data (i.e., data generated by simulating the model equations) illustrating the performance and robustness of the algorithm. In \S \ref{sec:real_data} we apply the algorithm to the identification of parameters in a model for the monopolar growth of fission yeast cells where the observations consist of experimental data of cells migrating in vitro.  Finally in \S \ref{sec:conclusion} we summarise our findings. make some conclusions and state possible directions for future work. Further details of the experimental setup used to generate the experimental data are given in Appendix \ref{Ap:bio_methods}.

\section{Problem formulation}\label{sec:setup}
Throughout $\G$ denotes a closed smooth oriented $d-1$ dimensional hypersurface in $\Reals^d$,  $d=2,3$, with  outward pointing unit normal $\normal$.  Given a function $\eta$ defined in a neighbourhood of $\G$, the tangential or surface gradient of $\eta$ denoted by $\nabla_\G$ is defined as
\begin{equation}
\nabla_\G\eta:=\nabla\eta-\nabla\eta\cdot\normal\normal,
\end{equation}
where $\nabla$ denotes the Cartesian gradient in $\Reals^d$. The Laplace-Beltrami operator $\lap_\G$ is defined as the tangential divergence of the tangential gradient, i.e.,
\begin{equation}
\lap_\G\eta:=\nabla_\G\cdot\left(\nabla_\G\eta\right).
\end{equation}
The mean curvature $H$ of $\G$ with respect to the normal $\normal$ is defined as
\begin{equation}
H:=\nabla_\G\cdot\normal.
\end{equation}
Note that by this definition the mean curvature is the sum of the principal curvatures and differs from the more common definition by a factor $\frac{1}{d-1}$. Note also that our sign convention is such that  the unit sphere has positive mean curvature if $\normal$ is the unit outer normal.

 The general evolution law we consider is of the form (for details of the modelling see \citet{venk11chemotaxis}) 
\beq\label{eqn:evolution_law}
\begin{split}
\vec V=&\Bigg(-\sigma H+k_b \left(\lap_\Gt{H}+H\lv\nabla_\Gt\normal\rv^2-\frac{1}{2}H^3\right)\\
&+g(\vec a)+\lambda(\mbox{Vol}(t)-\mbox{Vol}(0))\Bigg)\normal
\quad \text{on }\Gt,t\in(0,T],\\
\G(0)=&\G^0,
\end{split}
\eeq
where $\G$ is the closed surface that represents the cell membrane, $\vec V$ is the material velocity of $\G$, $\sigma$ is the surface tension,  $k_b$ is the bending rigidity, $\vec a$ is a vector of surface  resident species that satisfy (\ref{eqn:RDS}) a surface reaction-diffusion system (RDS), $g(\cdot)$ is the RDS species dependent forcing, $\mbox{Vol}(s)$ denotes the volume enclosed by the surface  at time $s$  and $\lambda\geq0$ is a spatially uniform penalisation term that, for strictly positive $\lambda$, accounts (weakly) for volume conservation. For the evolution of the surface resident species we shall consider an  RDS of the form
\beq
\label{eqn:RDS}
\begin{split}
\mdt_{\vec V}\vec a+\vec a\nabla_{\Gt}\cdot\vec V-\vec D\lap_{\Gt}\vec a&=\vec f(\vec a) \quad \text{on }\Gt,t\in(0,T],\\
\vec a(\cdot,0)&=\vec a^0(\cdot)\quad \text{on }\G(0),
\end{split}
\eeq
where $\vec{a} = (a_1, \dots, a_{n_a})^T$, $n_a$ is the number of chemical species involved, $a_i$ denotes the density of the $i$th chemical species, $\vec{V}$ is the material velocity of the surface, 
\begin{equation}\label{eqn:material_derivative}
\mdt_{\vec V}{\vec{a}}:=\pdt{\vec{a}}+\vec{V}\cdot\nabla{\vec{a}},
\end{equation}
is the material derivative with respect to the velocity $\vec V$, $\vec{D}$ is a diagonal matrix of positive diffusion coefficients and $\vec{f}(\vec{a})$ is the nonlinear reaction.

The aim of this work is to propose a method for the identification of parameters $\vec c\in\Reals^{n_p}$ such that the solution $(\Gt,\vec a(\vec x,t)), t\in(0,T],\vec x\in\Gt$ to (\ref{eqn:evolution_law}) and (\ref{eqn:RDS}) is ``close'' to some observed data. For example the parameters we wish to identify could correspond to the surface tension or bending rigidity $\sigma,k_b$ appearing in (\ref{eqn:evolution_law}), the diffusion coefficients $D_i$ appearing in (\ref{eqn:RDS}) or the forcing $g(\vec a)$ and the reaction $\vec f(\vec a)$ could be parameterised. 

In order for (\ref{eqn:evolution_law}) and (\ref{eqn:RDS}) to be well posed, it is natural assume the parameters take values in some admissible set $\mmm U_{ad}$. For example we could impose point wise bounds on the parameters, i.e, an admissible set of the form
\beq
\mmm U_{ad}:=\lbrace \vec c\in\Reals^{n_p}:\text{ for }i=1,\dotsc,n_p,\ \tilde c_i\leq c_i \leq \bar c_i\rbrace.
\eeq

The forward problem we consider in this study can now be stated as follows:
\begin{Pbm}[Forward problem]\label{pbm:forward}
Given parameters $\vec c\in\mmm U_{ad}$ and initial data  $(\G^0,\vec a^0)$  find $(\Gt,\vec a(\vec x,t))$ satisfying (\ref{eqn:evolution_law}) and (\ref{eqn:RDS}) for all $t\in(0,T],\vec x\in\Gt$.
\end{Pbm}

We now formulate the identification problem we shall consider in this study. We assume we have $n_s\geq 1$ observations of the data at times $0<t^1<..<t^{n_s}$. 
For a fixed $t^i$ we denote by $\left(\Gc(t^i),\hat{\vec a}(\vec x,t^i)\right),\vec x\in\Gc(t^i)$ the  associated observations. We denote by $(\Gt,\vec a(\vec x,t)), t\in(0,T],\vec x\in\Gt$,  the solution of the model  equations (\ref{eqn:evolution_law}) and (\ref{eqn:RDS}) with $T=t^{n_s}$.  Note we use the observation at $t=t^0$ to define the initial data for the model equations (\ref{eqn:evolution_law}) and (\ref{eqn:RDS}).
A key step in formulating the identification problem is to construct a suitable objective functional that measures the closeness of the solution to Problem \ref{pbm:forward} to the observations.   We propose two different formulations of the objective functional.
 
Firstly, we propose a sharp interface formulation of the objective functional of the following form 
\beq\label{eqn:obj_sharp}
\begin{split}
\mmm J_{sharp}&(\vec c)=
\frac{1}{2}\sum_{i=1}^{n_s}w_i\left(\int_{\G(t^i)} \lv d_{\Gc(t^i)}(\vec x)\rv \diff\vec x+\int_{\Gc(t^i)}\lv d_{\G(t^i)}(\hat{\vec x})\rv\diff\hat{\vec x}\right)^2\\
+w_{i+n_s}&\left(\int_{\G(t^i)}\lv\vec a(\vec x,t^i)-\hat{\vec{a}}(\hat{\vec x}_{cp}(\vec x,t^i),t^i)\rv\diff{\vec x}+\int_{\Gc(t^i)}\lv\hat{\vec a}(\hat{\vec x},t^i)-\vec{a}(\vec x_{cp}(\hat{\vec x},t^i),t^i)\rv\diff\hat{\vec x}\right)^2,
\end{split}
\eeq
where the $w_i\in\Reals^+$ are weights which may be tuned depending on the problem, $ d_\Gt(\cdot)$ and $ d_\Gct(\cdot)$ denote the (signed) distance functions to $\Gt$ and $\Gct$ respectively and
where for a given point $\vec p\in\Reals^{d}$, $\vec x_{cp}(\vec p,t)$ and  $\hat{\vec x}_{cp}(\vec p,t)$ denote the closest points on  $\Gt$ and  $\Gct$ respectively. 

Secondly, we consider a phase field formulation of the objective functional. Let $\O(t^i)\subset\Reals^{d}, i=1,\dots,n_s$ be such that for $i=1,\dots,n_s,$ 
\beq
\G(t^i)=\lbrace \vec x\in\O(t^i)|d_{\G(t^i)}(\vec x)=0\rbrace\quad\text{and}\quad\Gc(t^i)=\lbrace \vec x\in\O(t^i)|d_{\Gc(t^i)}(\vec x)=0\rbrace.
\eeq
Let $\phi^i_\epsilon,\hat{\phi}^i_\epsilon:\O(t^i)\to\Reals, i=1,\dots,n_s$ be such that for $i=1,\dots,n_s$  and $\vec x\in\O(t^i)$,
\beq\label{eqn:phi}
\phi^i_\epsilon(\vec x)=
\begin{cases}
1\quad&\text{if }d_{\G(t^i)}(\vec x)>\epsilon\\
\sin\left(\frac{\pi d_{\G(t^i)}(\vec x)}{2\epsilon}\right)\quad&\text{if }\lv d_{\G(t^i)}(\vec x)\rv<\epsilon\\
-1\quad&\text{if }d_{\G(t^i)}(\vec x)<-\epsilon
\end{cases}
\eeq
and
\beq\label{eqn:hat_phi}
\hat{\phi}^i_\epsilon(\vec x)=
\begin{cases}
1\quad&\text{if }d_{\Gc(t^i)}(\vec x)>\epsilon\\
\sin\left(\frac{\pi d_{\Gc(t^i)}(\vec x)}{2\epsilon}\right)\quad&\text{if }\lv d_{\Gc(t^i)}(\vec x)\rv<\epsilon\\
-1\quad&\text{if }d_{\Gc(t^i)}(\vec x)<-\epsilon,
\end{cases}
\eeq
where $\epsilon$ is a small positive parameter. Also let $\vec a^i_\epsilon,\hat{ \vec a}^i_\epsilon:\O(t^i)\to\Reals^m, i=1,\dots,n_s$ be such that for $i=1,\dots,n_s$ and $\vec x\in\O(t^i)$,
\beq\label{eqn:hat_a}
\vec a^i_\epsilon(\vec x)=
\begin{cases}
\cos\left(\frac{\pi d_{\G(t^i)}(\vec x)}{2\epsilon}\right)E_\G(t^i)[\vec a^i](\vec x,t^i)&\text{if }\lv d_{\G(t^i)}(\vec x)\rv<\epsilon\\
0&\text{otherwise},
\end{cases}
\eeq
and
\beq
\hat{\vec a}^i_\epsilon(\vec x)=
\begin{cases}
\cos\left(\frac{\pi d_{\Gc(t^i)}(\vec x)}{2\epsilon}\right)E_{\Gc(t^i)}[\hat{\vec a}](\vec x,t^i)&\text{if }\lv d_{\Gc(t^i)}(\vec x)\rv<\epsilon\\
0&\text{otherwise},
\end{cases}
\eeq
where $E_\Gt$ and $E_\Gct$ denote the constant normal extension operators to $\Gt$ and $\Gct$ respectively.
The phase field formulation of the objective functional we shall consider is given by
 \beq\label{eqn:obj_pf}
\begin{split}
\mmm J_{pf,\epsilon}(\vec c)=\frac{1}{2}\sum_{i=1}^{n_s}\left(w_i\ltwon{\hat{\phi}_\epsilon-\phi_\epsilon}{(\O(t^i))}^2+w_{i+n_s}\ltwon{\hat{\vec a}_\epsilon-\vec a_\epsilon}{(\O(t^i))}^2\right).
\end{split}
\eeq

The sharp interface formulation of the objective functional has the advantage that it is computed on the surfaces $\G$ and $\Gc$,  i.e., surfaces in $\Reals^{d-1}$, moreover, in \S \ref{subsec:objective_discretisation} we propose a discretisation of the sharp interface objective functional that is applicable to point cloud data which is important in applications. However, the sharp interface  formulation of the objective functional is not differentiable and the optimisation method we employ for the solution of the identification problem can only be shown to converge for smooth $(C^2)$ objective functionals. The phase field formulation of the objective functional, which is uniquely defined for $\epsilon$ sufficiently small, is evaluated on a bulk domain $\O\subset\Reals^{d}$ which adds to the computational cost of the method but it is smooth and there is therefore scope for the analysis of the identification method we propose in this study.  

 The identification problem  we shall consider may now be stated as follows:
 \begin{Pbm}[Identification problem]\label{pbm:id}
Given observations  $\left(\Gc(t^i),\hat{\vec a}(\vec x,t^i)\right),\vec x\in\Gc(t^i), i=0,\dots,n_s,$  find parameters $\vec c^*\in\mmm U_{ad}$, such that with $\left(\G^0,\vec a^0\right)$ defined by $\left(\Gc(t^0),\hat{\vec a}(\vec x,t^0)\right),\vec x\in\Gc(t^0)$ and   $(\Gt,\vec a(\vec x,t))$, $\vec x\in\Gt,t\in(0,T]$ solutions of Problem \ref{pbm:forward}, $\vec c^*$ solves the minimisation problem
\beq
\min_{\vec c\in\mmm U_{ad}}\mmm J(\vec c),\quad \text{with }\mmm J(\cdot)\text{ given by (\ref{eqn:obj_sharp}) or (\ref{eqn:obj_pf}).}
\eeq
\end{Pbm}

\section{Optimisation method}\label{sec:optimisation}
We propose a gradient-based optimisation algorithm for the solution of the optimisation problem, which requires the evaluation of the gradient of the objective functional.  Adjoint-based approaches are one commonly used method  to evaluate the gradient \citep{hinze2009optimization}. One advantage of such an approach is that  only one solution of the so called adjoint system (per optimisation iteration) is required for each evaluation of the gradient, independent of the number of parameters to be estimated. However, this approach necessitates the derivation of the optimality system, for which there is as yet no adequate theory in the present setting. Moreover, for free boundary problems one expects the solution to the forward problem, i.e., the geometry and the concentrations, to enter the adjoint equation which is posed backwards in time \citep{hogea2008image}, this can become computationally prohibitive  in terms of memory requirements especially in $3d$ \citep{hausser2010influence,hausser2012control}.
The gradient may be approximated with finite differences,  for each optimisation iteration this requires multiple forward problem solves equal to the number of parameters to be estimated. However,  the numerical results reported in  \citet{hogea2008image} suggest that the adjoint based method they consider is comparable in terms of CPU time to a finite difference gradient based method for the case that 3-5 parameters are to be estimated. Moreover, as each forward solve required for the finite difference approximation of the gradient is independent they may be carried out in parallel. An alternative to gradient based optimisation methods is the Bayesian approach \citep{stuart2010inverse}. This has been applied to parameter identification in reaction diffusion systems and has advantages in terms of obtaining confidence intervals for the estimated parameters \citep{dewar2010parameter}. Our focus in this study is to apply the efficient and robust solver we have developed for the forward problem, Problem \ref{pbm:forward}, to investigate the identification problem and hence we consider the use of a finite difference based gradient approximation in a gradient based optimisation method deferring to future studies the adjoint and Bayesian approaches. Specifically, we investigate the use of the Levenberg-Marquardt method for the solution of the optimisation problem, which is a widely used method that exploits the least squares structure of the problem and is thought to be robust to problems with large residuals  \citep{kelley1999iterative}.

For the general theory of parameter identification, optimisation and related inverse problems and optimal control problems we refer,  for example, to \citep{nocedal1999numerical, isakov1998inverse, troltzsch2010optimal}
and for theoretical results on parameter identification for semilinear evolution equations we refer to \citep{ackleh1998parameter}.

\subsection{The Levenberg-Marquardt algorithm}\label{sec:LM}
We seek to employ the Levenberg-Marquardt (LM) algorithm \citep{marquardt1963algorithm} for the solution of Problem \ref{pbm:id}, a standard method for nonlinear least squares problems which has been applied in parameter identification for partial differential equations \citep{burger2004levenberg, iglesias2011level}. We now present the LM algorithm in this context, for further details of the applications and analysis of the algorithm we refer, for example to  \citet{more1978levenberg}.

We define $\vec \chi\in\Reals^{2n_s}$ such that the objective functional $\mmm J=\frac{1}{2}\sum^{2n_s}_{i=1}\chi_i^2=\frac{1}{2}\vec \chi^T\vec \chi$, with the components $\chi_i$ given by for $ i=1,\dots,n_s$
\beq\label{eqn:chi_sharp}
\begin{split}
\chi_i(\vec c)=&
w_i^{1/2}\left(\int_{\G(t^i)} \lv d_{\Gc(t^i)}(\vec x)\rv\diff\vec x+\int_{\Gc(t^i)}\lv d_{\G(t^i)}(\hat{\vec x})\rv\diff\hat{\vec x}\right)\\
\chi_{i+n_s}(\vec c)=&
w_{i+n_s}^{1/2}\Bigg(\int_{\G(t^i)}\lv\vec a(\vec x,t^i)-\hat{\vec{a}}(\hat{\vec x}_{cp}(\vec x,t^i),t^i)\rv\diff{\vec x}+\int_{\Gc(t^i)}\lv\hat{\vec a}(\hat{\vec x},t^i)-\vec{a}(\vec x_{cp}(\hat{\vec x},t^i),t^i)\rv\diff\hat{\vec x}\Bigg)
\end{split}
\eeq
in the sharp interface objective functional case (\ref{eqn:obj_sharp}) and
\beq\label{eqn:chi_pf}
\chi_i(\vec c)=
\begin{cases}
w_i^{1/2}\ltwon{\hat{\phi}_\epsilon-\phi_\epsilon}{(\O(t^i))}\quad i=1,\dots,n_s\\
w_i^{1/2}\ltwon{\hat{\vec a}_\epsilon-\vec a_\epsilon}{(\O(t^{i-n_s}))}\quad i=n_s+1,\dots,2n_s,
\end{cases}
\eeq
in the phase field objective functional case (\ref{eqn:obj_pf}).

We denote by $\vec J=\vec J(\vec c)\in\Reals^{2n_s\times n_p}$ the Jacobian matrix of $\vec \chi$ with components 
\beq\label{eqn:Jac}
J_{ij}=\frac{\partial \chi_i}{\partial c_j}.
\eeq 
Similar to the Gauss-Newton method, the LM method is based on considering  a linearisation of $\vec \chi$ and approximating the Hessian of the objective functional $\mmm J$ with $\vec  J^T\vec J$.  For each iteration of the algorithm, given an initial guess for the parameters $\vec c$ an update $\vec \delta\in\Reals^{n_p}$ is computed by solving 
$$
(\vec J^T\vec J +\mu \vec I)\vec \delta 
= 
-\vec J^T\vec\chi,
$$
where $\vec J=\vec J(\vec c),\vec \chi =\vec\chi(\vec c)$. The update is then
$$
\vec c = \vec c +\vec \delta.
$$
The damping $\mu\in\Reals^+$ differentiates LM from Gauss-Newton. For large values of $\mu$ LM resembles steepest descent while for smaller values of $\mu$ LM resembles Gauss-Newton. In practice $\mu$ is chosen adaptively with the value of $\mu$ decreased whilst each iteration of the algorithm results in a decrease of the error.
The termination conditions for the algorithm are given in terms of
\begin{itemize}
\item The magnitude of the objective functional  $\vert \mmm J\vert$.
\item The magnitude of the gradient of the linearised objective functional.
\item The relative change in the update $ \vec\delta_{LM}$.
\item A maximum number of iterations.
\end{itemize}
A necessary condition for the approximation of the Hessian $\vec J^T\vec J$ to be invertible is that we have at least as many observations as we have parameters we wish to identify. Thus in the present setting we require $2n_s\geq n_p$.

\section{Discretisation and implementation}\label{sec:implementation}
For the problems we have in mind, it is generally not possible to determine the components of the Jacobian matrix $\vec J$, c.f., (\ref{eqn:Jac}), analytically. We therefore employ finite differences to approximate $\vec J$ \citep[Chap. 8.6]{monahan2011numerical}. This procedure necessitates multiple evaluations of the functional $\vec \chi$. We note that the evaluations of the the functional $\vec \chi$ needed to compute a finite difference approximation to the Jacobian matrix $\vec J$ are independent of each other and may therefore be carried out in parallel (initial investigations in this direction show that a  parallel implementation gives  significant computational savings for the Jacobian evaluation  stage \citep{Liuidparallel}). 

As analytical solutions to Problem \ref{pbm:forward} are generally unavailable we must approximate the solution to Problem \ref{pbm:forward}. We employ a robust and efficient method based on the surface finite element method \citep{dziuk2008computational,barrett2008parametric,doi:10.1137/110828642,dziuk2007finite,dziuk2010l2,lubich2013backward,EllSty12} for this purpose. The numerical method we employ for the solution of the forward problem and the applications of this method to the study of cell motility are described in detail elsewhere \citep{venk11chemotaxis}. 

The LM algorithm is implemented with the help of the levmar library  \citep{lourakis2004levmar}. The library provides routines for the LM algorithm with finite difference evaluation of the Jacobian and incorporates box constraints (i.e., pointwise bounds on the parameters).

We stress that the forward problem is solved with a parametric finite element method even when the phase field formulation of the objective functional (\ref{eqn:obj_pf}) is used. Thus for each time at which we want to compute the objective functional (\ref{eqn:obj_pf}), we construct a phase field representation of the computed and target surfaces together with the corresponding (weighted by distance constant normal) extensions  of the target concentrations and computed concentrations.

\subsection{Discretisation of the objective functional}\label{subsec:objective_discretisation}
The numerical method we employ for the approximation of the solution to Problem \ref{pbm:forward} consists in approximating the surface $\G$ with a triangulated surface. For simplicity in this study, we focus on the case where the surface $\G$ is approximated by a piecewise linear triangulated surface $\G_h$. We denote by $h$ the mesh-size of the triangulated surface $\G_h$.

Our methodology is formulated assuming a continuous description of  the observed data, however, in practical applications it is often the case that the observations are not given as a series of surfaces with associated concentrations, rather at each time the observed data is a point cloud consisting of points that lie on the cell membrane together with associated concentrations at these points.

In the case that we consider the sharp interface objective functional $\mmm J_{sharp}$ we compare the positions and associated concentrations at the vertices of the computed surfaces directly with the point cloud observations. To this end we introduce the Haussdorff distance, given a set of points $\vec P:=\lbrace\vec p_i:i=1,\dotsc,n_p\rbrace$ with $\vec p_i\in\Reals^{d}$ for each $i$ and a point $\vec x\in\Reals^{d}$, we define the Hausdorff distance $d^H_{\vec P}(\vec x)$  between the point $\vec x$ and the set of points $\vec P$,  by
\beq
d^H_{\vec P}(\vec x):=\inf_i\lv \vec x-\vec p_i\rv.
\eeq
Let $N_{\G_h}$ and $N_{\Gc_{\hat h}}$ denote the number of vertices of the triangulated surface and the number of points in the point cloud observations respectively  and let $d^H_{\G_h}$ and $d^H_{\Gc_{\hat h}}$ denote the Hausdorff distance to the set of vertices of the triangulated surface  and the point cloud of observations respectively. We propose the following discretisation of the objective functional (\ref{eqn:obj_sharp}): 
\beq\label{eqn:disc_obj_sharp}
\begin{split}
\mmm J_{h,sharp}&(\vec c)=
\sum_{i=1}^{n_s}w_i\frac{1}{2}\left(\sum_{j=1}^{N_{\G_h}}  d^H_{\Gc_{\hat h}(t^i)}(\vec x_j)/N_{\G_h}+\sum_{k=1}^{N_{\Gc_{\hat h}}} d^H_{\G_h(t^i)}(\hat{\vec x}_k)/N_{\Gc_{\hat h}}\right)^2\\
+w_{i+n_s}\frac{1}{2}&\left(\sum_{j=1}^{N_{\G_h}}\lv \vec a(\vec x_j,t^i)-\hat{\vec{a}}(\hat{\vec x}^H_{cp}(\vec x_j,t^i),t^i)\rv/N_{\G_h}+\sum_{k=1}^{N_{\Gc_{\hat h}}}\lv\hat{\vec a}(\hat{\vec x}_k,t^i)-\vec{a}(\vec x^H_{cp}(\hat{\vec x}_k,t^i),t^i)\rv/N_{\Gc_{\hat h}}\right)^2,
\end{split}
\eeq
where the closest point operators ${\vec x}^H_{cp}$ and $\hat{\vec x}^H_{cp}$ are the closest points in the set of vertices of the triangulated surface and the point cloud observations respectively. We observe that as the discretisation parameters $h$ and $\hat h$ tend to zero, i.e., as the number of distinct points in the point cloud of observations and the number of vertices in the triangulated surface both tend to infinity, we recover the continuous objective functional $\mmm J_{sharp}$ (c.f., (\ref{eqn:obj_sharp})), however, in a practical setting the observations are likely to be (potentially noisy) data generated in experiments with a fixed resolution and we may only have control of the parameter $h$ and not the discretisation parameter associated with the observations $\hat h$.

For the phase field objective functional (\ref{eqn:obj_pf}), if we are given observations in the form of point cloud data , as we need a signed distance function we must construct a surface from the target data. The construction of triangulated surfaces from point clouds is not the focus of this work and therefore we assume that if the data is given in the form of a point cloud then we have suitable information on connectivity such that the points in the cloud represent  the vertices of  a triangulated surface. Given a pair of  triangulated surfaces (one for the observations and one from the computations) we  define a rectangular domain $\O$ such that the minimum distance between any point on the triangulated surfaces and the boundary of the rectangular domain is greater than $2\epsilon$. We then triangulate $\O$ such that $h_\O<\frac{\epsilon}{4}$, where $h_\O$ is the maximal mesh size of all the elements of  this (bulk) triangulation. We then   approximate the phase field representation of the position and the concentrations (c.f., (\ref{eqn:phi})---(\ref{eqn:hat_a})) with piecewise linear $C^0$ functions (i.e., $\mathbb{P}^1$ finite element functions) over this bulk triangulation.


\section{Numerical experiments with artificial data}\label{sec:examples}
We focus on an example that arises in the modelling of the motility of Keratocyte fragments \citep{venk11chemotaxis}. We consider   Problem  \ref{pbm:forward} with a 2-component RDS. For the forcing function $g$ appearing in the evolution law (\ref{eqn:evolution_law}) we assume protrusion or retraction proportional to the RDS species, i.e., with $\vec k_p\in\Reals^2$,
\beq\label{eqn:prop_forcing}
g(\vec a)=\vec k_p\cdot \vec a.
\eeq
 We choose the reaction kinetics appearing in \ref{eqn:RDS} to be given by
\beq\label{eqn:schnak}
f_1(\vec a)=\gamma(k_1-a_1+a_1^2a_2)\quad\mbox{and}\quad f_2(\vec a)=\gamma(k_2-a_1^2a_2).
\eeq
The reaction kinetics (\ref{eqn:schnak}) are known as the activator-depleted substrate model \citep{prigo68}. They have been used widely in the modelling of biological pattern formation phenomena and there has been  progress in the analysis of RDSs on fixed and growing domains equipped with the kinetics (\ref{eqn:schnak}) at the continuous level \citep{kol2008ssr,Venkataraman2010global,venkataramanthesis} and in the analysis of numerical methods to approximate the equations \citep{garvie2010efficient,venk10fem,venkataraman2013adaptive}.

In all the numerical examples in this section we take either the unit circle or unit sphere as the initial position of the curve or surface at $t=0$. 
We select the volume conservation penalisation term $\lambda$,  appearing in (\ref{eqn:RDS}), to be equal to 1.
The initial conditions for the surface RDS were taken to be the spatially homogeneous steady state value for the activator species, i.e.,  $a_1(\vec x,0)=1$ with a perturbation of the form $a_2(\vec x,0)=0.9+0.001\max(0,-x_1)$ introduced to the substrate concentrations.

For the LM algorithm, for all the results in this Section, we used the parameter values given in Table \ref{tab:LM_parameter_values}. We employed box constraints on all of the parameters such that they were each constrained to be positive or negative depending on the sign of the true value and their magnitude  was constrained to be less than three times the magnitude of the true value. 

In all the simulations apart from those in \S \ref{subsec:sensitivity}, we set the weights $w_i$, appearing in (\ref{eqn:chi_sharp}) and (\ref{eqn:chi_pf}), to be $1$ for each $i$ and  for the formulation of the phase field  objective functional we took $\epsilon$ (the parameter that governs the interface width) to be equal to 0.5.
To deal with parameters with differing orders of magnitude, for the finite difference step size we rescaled each of the parameters to be identified such that the true value was 1. Thus a step size of $1\times10^{-2}$ for a parameter whose true value was equal to $10$ or $0.1$ would correspond to an actual step size of $0.1$ or $1\times10^{-3}$ respectively in the finite difference Jacobian approximation step.
\begin{table}
\centering
\begin{tabular}{cccc}
\toprule
&\multicolumn{3}{c}{Stopping criteria}\\
$h_{LM}$&Gradient: $\left\| \vec J^T \vec \chi\right\|_{\infty}$&Update: $\|\vec \delta\|_2$&Error: $\|\vec \chi\|_2$\\
\midrule
$5\times10^{-3}$&$1\times10^{-6}$&$1\times10^{-6}$&$1\times10^{-6}$\\
\bottomrule
\end{tabular}\\
\caption[]{Parameter values for the Levenberg-Marquardt algorithm for all the experiments reported on in \S  \ref{sec:examples}. The notation corresponds to the notation introduced in \S \ref{sec:LM}.}\label{tab:LM_parameter_values}
\end{table}

In the majority of the examples we report on, the algorithm terminated due to the estimate of the error (i.e., the objective functional) being below the given tolerance, i.e.,  $\|\vec \chi\|_2<1\times10^{-6}$, we note that in the case of experimental data we would not expect this to be the case as the model error is unlikely to be negligible. In the  example corresponding to the estimates obtained using the phase field objective functional stated in the second part of Table \ref{tab:surface_identification} and the majority of examples computed with noisy data (Table \ref{tab:noisy_curve_identification}) the algorithm terminated due to the estimate of the gradient being below the given tolerance, i.e.,  $\left\| \vec J^T \vec \chi\right\|_{\infty}<1\times10^{-6}$.

\subsection{Numerical experiments for curves}\label{subsec:curves} 
We start with an experiment where the surface  $\G\subset\Reals^2$ is a closed curve. We generate the target data by approximating the solution to Problem \ref{pbm:forward} with a forcing term of the form (\ref{eqn:prop_forcing}) and reaction kinetics (\ref{eqn:schnak}). We selected  parameter values for the RDS and surface evolution law stated in Table \ref{tab:curve_parameter_values}. We approximated the solution to Problem \ref{pbm:forward} using linear finite elements on a mesh with 128 degrees of freedom and selected a timestep of $10^{-2}$.
\begin{table}
\centering
\begin{tabular}{cccccccccc}
\toprule
$T$&$\sigma$&$k_b$&$D_1$&$D_2$&$(k_p)_1$&$(k_p)_2$&$\gamma$&$k_1$&$k_2$\\
\midrule
10&$5\times10^{-3}$&$1\times10^{-2}$&$1$&$100$&$-1\times10^{-2}$&$5\times10^{-2}$&$20$&$0.1$&$0.9$\\
\bottomrule
\end{tabular}\\
\caption[]{Parameter values used to generate the target data for numerical experiments on curves.}\label{tab:curve_parameter_values}
\end{table}

For the identification problem we attempted to recover the coupling term $(k_p)_2$ appearing in the forcing function $g$ that appears in the evolution law (\ref{eqn:evolution_law}),  the parameter $\gamma$ that appears in the reaction kinetics (\ref{eqn:schnak}) and the bending rigidity of the cell membrane $k_b$  that appears in the evolution law (\ref{eqn:evolution_law}). A possible physical interpretation of the parameter values is that $(k_p)_2$ corresponds to the protrusive force generated by actin filaments, $\gamma$ is the rate constant of the reactions taking place on or near the cell membrane that lead to polarisation and $k_b$ governs the magnitude of the force generated due to the resistance of the cell membrane to bending \citep{venk11chemotaxis}. We assumed the target data was observed at $t=0,1,2,\dots,10$ and thus had $10$ observations of the data to fit to with the initial observation used to define the initial data for the problem. For the phase field  objective functional we computed the objective functionals on rectangular domains such that the distance between the computed and target curves and the boundary of the domain was at least $2\epsilon$ and we used a triangulation of the rectangular domains with 4225 degrees of freedom. 

Table \ref{tab:curve_identification} shows the results of two experiments, one for each choice of the objective functional (\ref{eqn:disc_obj_sharp}) or (\ref{eqn:obj_pf}). The starting values for the parameters, used as an initial guess for the algorithm, together with the computed values are reported together with the relative percentage errors, where the relative percentage error is defined as 
\[
\lv\frac{\mbox{estimated value}-\mbox{true value}}{\mbox{true value}}\rv\times100.
\]
  The algorithm converged in 6 and 7 iterations for the sharp interface and phase field objective functionals respectively and took 142 and 197 seconds of CPU time for the sharp interface and phase field objective functionals respectively. The algorithm appears to perform best in both cases at  identifying the protrusive forcing strength $(k_p)_2$ and the reaction rate $\gamma$ with the relative error for the estimated value of the bending rigidity $k_b$ somewhat higher.
  
\begin{figure}[ht!]
\centering
\includegraphics[trim = 0mm 0mm 0mm 0mm, clip=true, width=\linewidth]{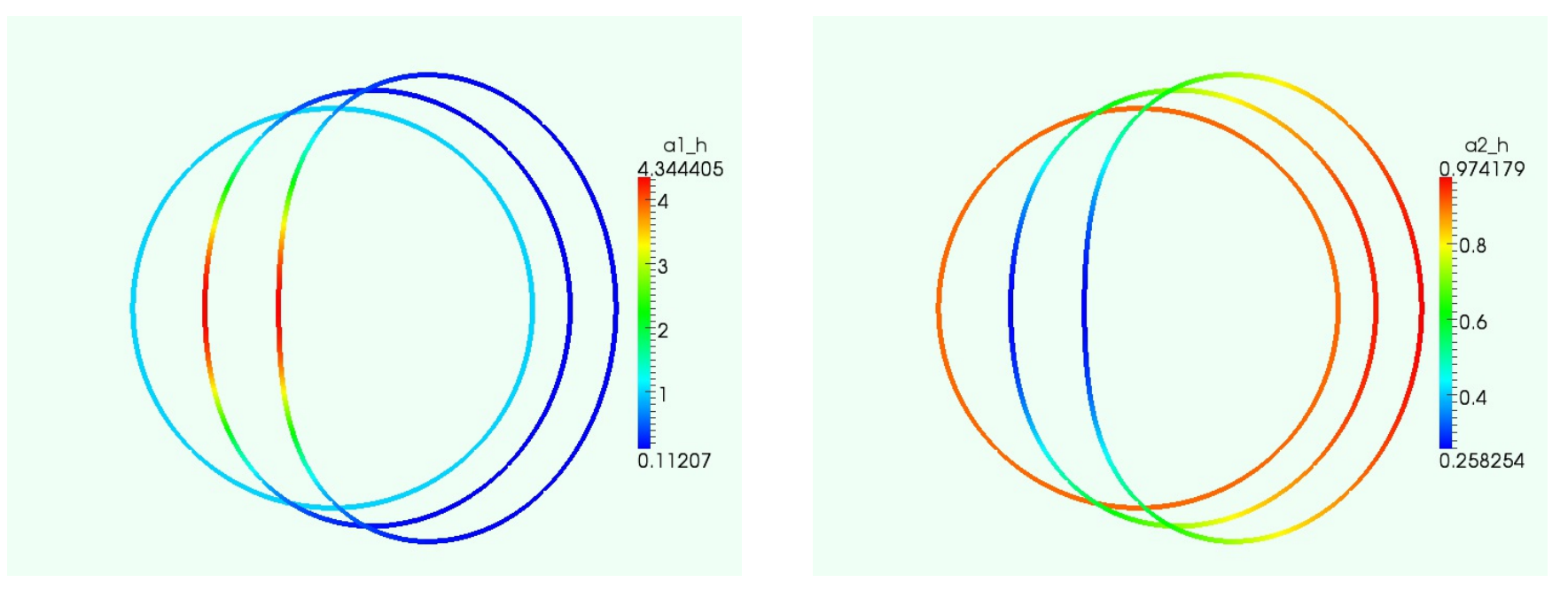}
\caption{(Colour online) Results of the simulation described in \S \ref{subsec:curves} with estimated parameters, obtained using the sharp interface objective functional, given in Table \ref{tab:curve_identification} at $t=0,5$ and $10$. Note the results are indistinguishable from those computed with the true parameter values (not reported). }\label{fig:sharp_curve_simulations}  
\end{figure}

\begin{figure}[ht!]
\centering
\includegraphics[trim = 0mm 0mm 0mm 0mm, clip=true, width=\linewidth]{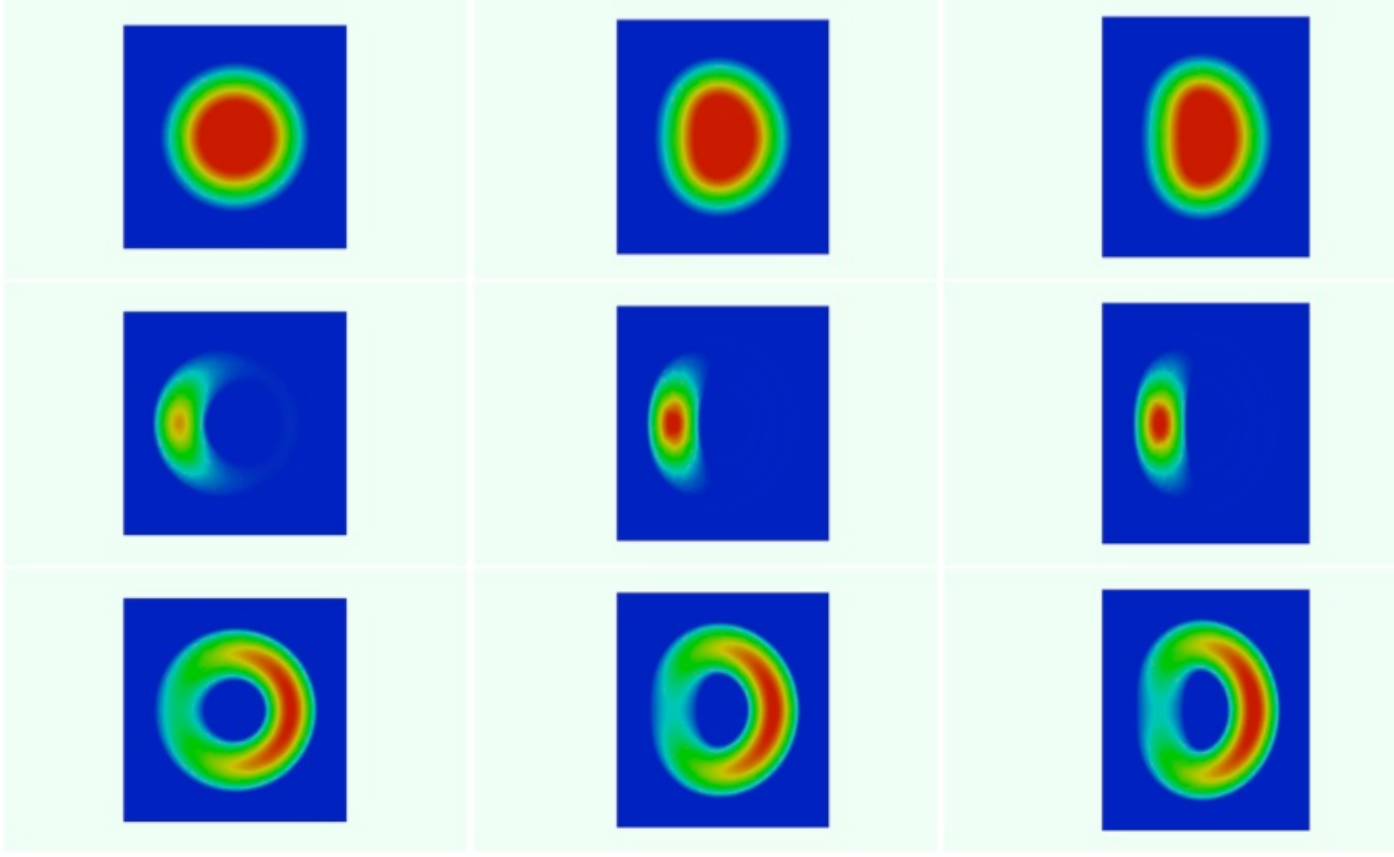}
\caption{(Colour online) Phase field representation of the simulation described in \S \ref{subsec:curves} with the estimated parameters given in Table \ref{tab:curve_identification} obtained using the phase field objective functional.  The position $\phi_\epsilon$ (c.f., (\ref{eqn:hat_phi})) (top row) and the RDS species $(a_1)_\epsilon$ (middle row) and $(a_2)_\epsilon$ (bottom row) are shown at $t=0,5,10$ reading from left to right. Note the results are indistinguishable from those computed with the true parameter values (not reported). }\label{fig:pf_curve_simulations}  
\end{figure}

  Figures \ref{fig:sharp_curve_simulations} and \ref{fig:pf_curve_simulations} show results of simulations with the estimated parameters using the sharp interface and phase field objective functionals respectively. The results are indistinguishable from those obtained using the true parameter values which we do not report on. In Figure \ref{fig:sharp_curve_simulations} we show the  RDS species concentrations on the triangulated surface $\G_h$ at $t=0,5$ and $10$ and in Figure \ref{fig:pf_curve_simulations} we show the  phase field formulation of the position of the surface and the RDS species concentrations  at $t=0,5$ and $10$. We note that although the problem is one in which both large deformations in the surface position and large changes in the concentrations of the RDS species are evident, the parameter identification algorithm successfully identifies parameters with small relative errors and for which the results of simulating the model equations are indistinguishable from those obtained using the true parameter values.
   
\begin{table}
\centering
\begin{tabular}{ccccc}
\toprule
&Parameter&True value&Starting value &Computed value\\
\midrule
\multirow{3}{0.25\textwidth}{Sharp interface objective functional: $\mmm J_{sharp}$}
&$(k_p)_2$&$5.00\times10^{-2}$&$3.75\times10^{-2}\ (25)$&$5.00\times10^{-2} \  (0.02)$\\
&$\gamma$&$20.00\times10^{0}$&$25.00\times10^{0} \ (25)$&$19.99\times10^{0} \  (0.04)$\\
&$k_b$&$1.00\times10^{-2}$&$1.15\times10^{-2}\ (15)$&$9.68\times10^{-3} \  (3.19)$\\
\\
\multirow{3}{0.25\textwidth}{Phase field objective functional: $\mmm J_{pf,\epsilon}$}
&$(k_p)_2$&$5.00\times10^{-2}$&$3.75\times10^{-2}\ (25)$&$5.00\times10^{-2} \  (0.05)$\\
&$\gamma$&$20.00\times10^{0}$&$25.00\times10^{0}\ (25)$&$19.99\times10^{0} \  (0.07)$\\
&$k_b$&$1.00\times10^{-2}$&$1.15\times10^{-2}\ (15)$&$9.88\times10^{-3} \  (1.22)$\\
\bottomrule
\end{tabular}\\
\caption[]{Results for an identification experiment for the forced motion of a curve, relative $\%$ error in parentheses.}\label{tab:curve_identification}
\end{table}             

\subsubsection{Noisy observations} 
In practical applications the observations are likely to include some noise either from measurement error or due to natural biological variability. To investigate the robustness of the proposed algorithm in identifying parameters from noisy data we investigated the performance of the algorithm when the observations were perturbed by noise. We denote by $\lbrace{\vec x}^i_k\rbrace_{k=1,\dotsc,N_{\Gc_h}}$  the set of vertices of the target curves at time $i$ and by  $\lbrace {\vec a}_k^i\rbrace_{k=1,\dotsc,N_{\Gc_h}}$ the associated concentrations. As the different RDS species and the spatial coordinates have very different scales, we consider perturbations of the following form, for $i=1,\dotsc,n_s,$ and for  $k=1,\dotsc,N_{\Gc_h},$
\beq\label{eqn:noise}
\begin{split}
(\hat{x}_j)_i^k&=({x}_j)_i^k+ \eta(\max_{l}(\hat{x}_j)_i^l-\min_{l}(\hat{x}_j)_i^l) \quad j=1,\dotsc, d,\\
(\hat{a}_j)_i^k&=({a}_j)_i^k+\eta (\max_{l}(\hat{a}_j)_i^l-\min_{l}(\hat{a}_j)_i^l)\quad j=1,\dotsc,n_a,
\end{split}
\eeq
where $d$ is such that $\G\subset\Reals^{d}$, $n_a$ denotes the number of RDS species and  $\eta$ is a random variable.  We consider the case where $\eta$ follows either  a normal distribution or a uniform distribution with mean zero and  standard deviation $k_n$.  We investigate the effect of varying the standard deviation $k_n$, i.e., the strength of the noise.

 Table \ref{tab:noisy_curve_identification} shows the results of the identification algorithm with observations perturbed by uniformly or normally distributed noise with varying standard deviation.  The starting values and true values of the parameters were as in Table \ref{tab:curve_identification}. We report the mean and standard deviation of the relative error for each parameter estimated after $100$ runs of the algorithm with noisy observations.  We observe that the results are similar for uniformly distributed and normally distributed noise with the same standard deviation. 
 
\begin{figure}[ht!]
\centering
\includegraphics[trim = 0mm 0mm 0mm 0mm, clip=true, width=\linewidth]{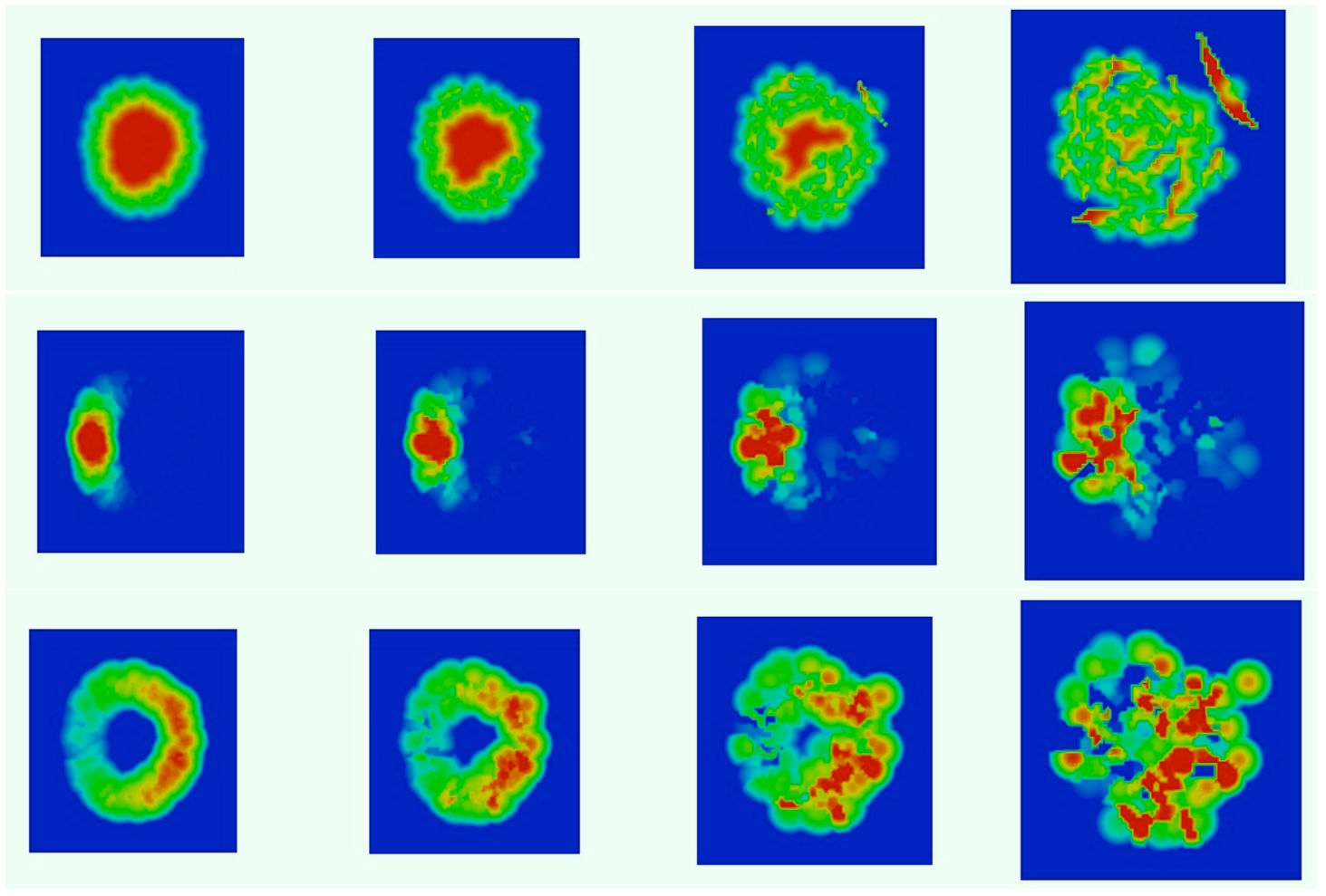}\\
\caption{(Colour online) Phase field representation the observations of the, position $\hat\phi_\epsilon$ (c.f., (\ref{eqn:hat_phi})) (top row) and the RDS species $(\hat a_1)_\epsilon$ (middle row) and $(\hat a_2)_\epsilon$ (bottom row) at $t=5$ for one simulation perturbed by normally distributed  noise. In each row reading from left to right the standard deviation of the noise  $k_n$ (c.f. (\ref{eqn:noise})) is $0.02, 0.05, 0.1$ and $0.2$. For comparison with the true values see middle column of Figure \ref{fig:pf_curve_simulations}. }\label{fig:noise_curve_pf_positon_simulations}  
\end{figure}

 The algorithm with the sharp interface objective functional appears to generate good estimates of the reaction rate $\gamma$ even with moderate levels of noise (up to standard deviations of $0.10$), for low to moderate noise the estimate of the protrusive forcing strength $(k_p)_2$  has a mean relative error  of less than $10\%$, while the estimate of the bending rigidity $k_b$ is not robust  even for noise with a small standard deviation and exhibits a large relative error (even greater than that of the starting value). On the other hand, the algorithm with the phase field objective functional generates estimates of the protrusive forcing strength and reaction rate that appear more sensitive to noise than the sharp interface case, however, the estimate of the bending rigidity appears more robust to noise with a relative error of around $10\%$ for small to moderate noise (standard deviations of $0.02-0.05$).  To illustrate the effect of the  perturbations by noise on the observations, in Figure \ref{fig:noise_curve_pf_positon_simulations} we report on the perturbed phase field representation of the observations at $t=5$  in one simulation with normally distributed noise with standard deviations of $0.02,0.05,0.1$ and $0.2$ respectively (the results are similar for uniformly distributed noise and are not reported on). We clearly observe that for noise with standard deviation of $0.1$ or greater that the observations are strongly perturbed from the true  observations which are indistinguishable from the middle column of Figure \ref{fig:pf_curve_simulations}.
 
In practice, it may be the case that experimentalists have knowledge on the noise inherent in given observations, for example the noise may increase with observations further in the future or due to measurement errors associated with certain snapshots. This could then be incorporated by tuning the weights such that larger relative weights were given to the observations that experimentalists had the most confidence in.

\begin{table}
\centering
\begin{tabular}{llcccc}
\toprule
&&Standard deviation&\multicolumn{3}{c}{Mean (and standard deviation) of the relative error}\\
&&$k_n$&$(k_p)_2$&$\gamma$&$k_b$\\
\midrule
\multirow{7}{0.1\textwidth}{Sharp interface}
&
\multirow{3}{0.12\textwidth}{Normally distributed noise}
&0.02&2.5037 (1.734)&0.1676 (0.132)&23.4624 (25.392)\\
&&0.05&7.0463 (5.515)&0.4813 (0.377)&34.5578 (24.835)\\
&&0.10&14.3527 (10.977)&    1.1055  (0.746) & 41.6302 (31.746)\\
&&0.20&   27.8036 (22.147) &  10.3681 (10.593)&   53.1654 (47.248)\\
\\
&
\multirow{3}{0.12\textwidth}{Uniformly distributed noise}
&0.02& 3.0980 (2.280)&    0.2398 (0.179)&   34.4362 (33.668)\\
&&0.05&  9.6986  (6.830)&   0.5784 (0.441)&  37.6576 (29.567)\\
&&0.10&    22.1223  (15.587)&  1.5397  (1.115)&  43.9082 (30.439)\\
&&0.20&   36.1964  (31.893)& 18.2845 (12.291)&   55.8188 (61.621)\\
\\
\multirow{7}{0.1\textwidth}{Phase field}
&
\multirow{3}{0.12\textwidth}{Normally distributed noise}
&0.02&4.9058 (1.606) &    3.6848 (0.773) &  11.2457 (9.277)\\
&&0.05&14.8669 (3.284) &   9.5646 (1.500) &   11.6231 (11.123) \\
&&0.10&30.7799 (8.179)&   17.3310  (3.644)&   18.9066 (14.866) \\
&&0.20&53.3621 (17.532) &  28.3924 (10.151)&   67.5399 (58.756) \\
\\
&
\multirow{3}{0.12\textwidth}{Uniformly distributed noise}
&0.02& 4.7467 (1.448) &    3.6671 (0.653) &  11.2149 (9.094)\\
&&0.05&  15.4467 (3.125)&    9.9262  (1.444)&  10.6822 (9.549)\\
&&0.10&   32.1184 (6.834)&   17.6269 (3.834)&  20.7039 (17.689)\\
&&0.20& 54.7630 (25.612)&  34.3233 (13.909)&   77.5613 (68.399)\\
\bottomrule
\end{tabular}\\\caption[]{Mean and standard deviation (in parentheses) of the relative error ($\%$)  of 100 runs of an identification experiment for the forced motion of a curve with normally and uniformly distributed noisy observations of varying standard deviation.}\label{tab:noisy_curve_identification}
\end{table}                    

\subsubsection{Sensitivity of the objective functional} \label{subsec:sensitivity} 
We now report on the results of experiments in which we plot the objective functional for a fixed set of observations generated using the true parameter values reported in Table \ref{tab:curve_identification}. We computed the value of the objective functional on a rectangular grid of 1028 points $\pm 25\%$ the true parameter value varying the parameters corresponding to the reaction rate $\gamma$ and the protrusive forcing strength $(k_p)_2$. 

\begin{figure}[hp]
\centering
\includegraphics[trim = 0mm 0mm 0mm 0mm, clip=true, width=0.8\linewidth]{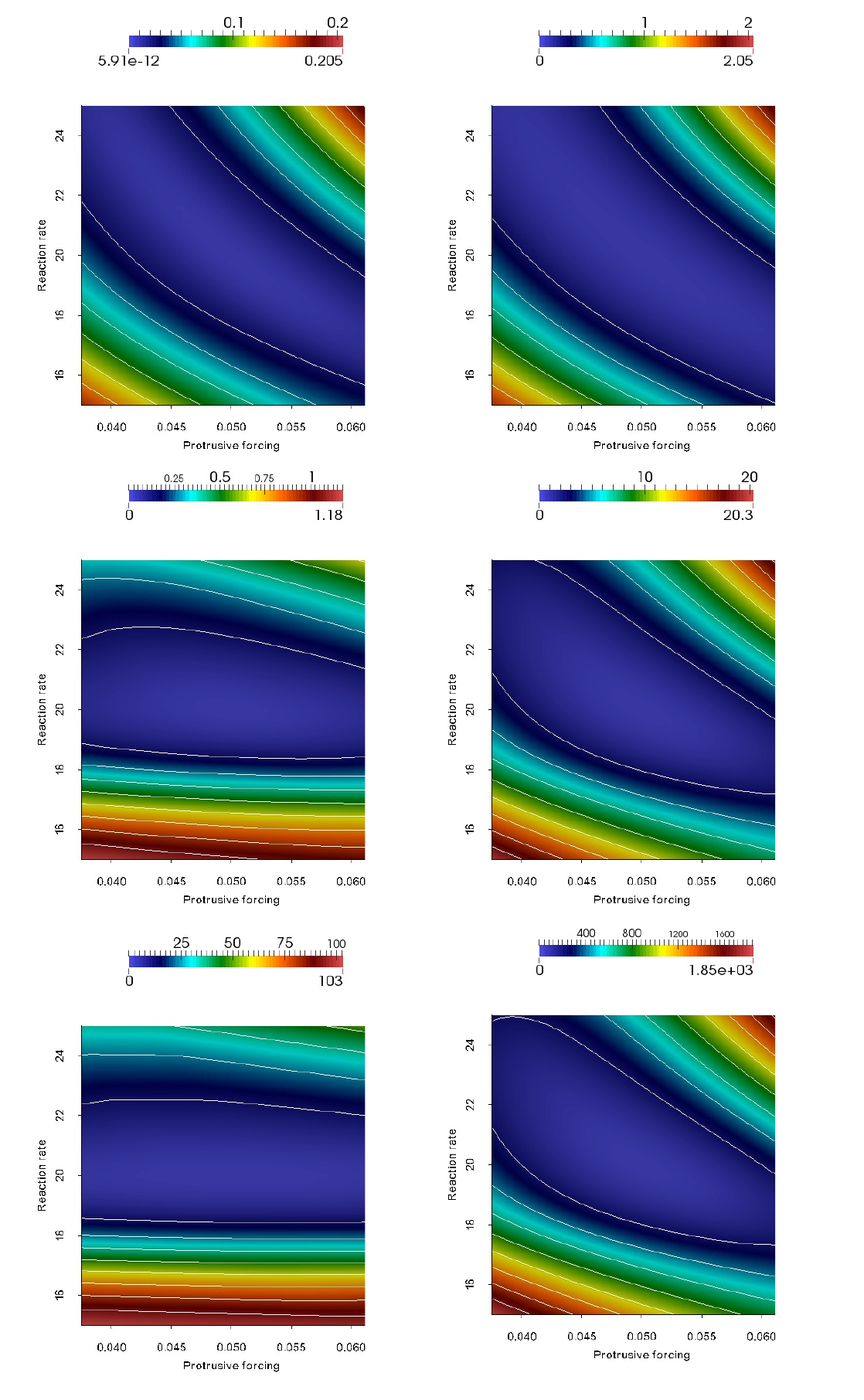}
\caption{(Colour online) Left: the sharp interface objective functional (\ref{eqn:disc_obj_sharp}) and right: the phase field objective functional (\ref{eqn:obj_pf}) as a function of the reaction rate $\gamma$ and the  protrusive forcing strength $(k_p)_2$.   The relative weight of the error due to the concentration $\alpha=0.01, 1, 100$, reading from top to bottom in each column (for details see text). The observations were generated using $\gamma=20$ and $(k_p)_2=0.05$. The lines indicate level curves of the objective functional.  }\label{fig:objective}
\end{figure}

We start by examining the sensitivity of the objective functional as we change the weights in (\ref{eqn:obj_sharp}) and (\ref{eqn:obj_pf}) such that the  contribution of the error due to position is changed relative to the error due to the concentrations. To do this we set $w_i=1$  and $w_{i+n_s}=\alpha$ for $i=1,\dots,n_s$ with $\alpha=0.01,1$ and $100$. Figure \ref{fig:objective} shows the dependence of the objective functional on the parameters for the sharp interface and phase field formulation.  When the contribution of the error due to concentration is small relative to the error due to position, i.e, $\alpha=0.01$, the sensitivity of the objective functional appears similar for the two different formulations (top row of Figure \ref{fig:objective}). For this choice of $\alpha$, the objective functionals appear relatively insensitive to a change in parameter values that corresponds to an increase in $\gamma$ combined with a reduction of $(k_p)_2$ as the level curves appear oriented along a line $y=-cx$, for some positive number $c$. In the sharp interface case (left column of Figure \ref{fig:objective}), as we increase the value of $\alpha$ the objective functional appears significantly more sensitive to changes in the reaction rate than the protrusive forcing strength with the level curves of the objective functional appearing almost horizontal. This is different to the experiments with the phase field objective functional where increasing the value of $\alpha$ only leads to slightly more horizontal level curves of the objective functional which appear to remain oriented along a line $y=-cx$, for some positive number $c$.

\begin{figure}[ht!]
\centering
\includegraphics[trim = 0mm 0mm 0mm 0mm, clip=true, width=\linewidth]{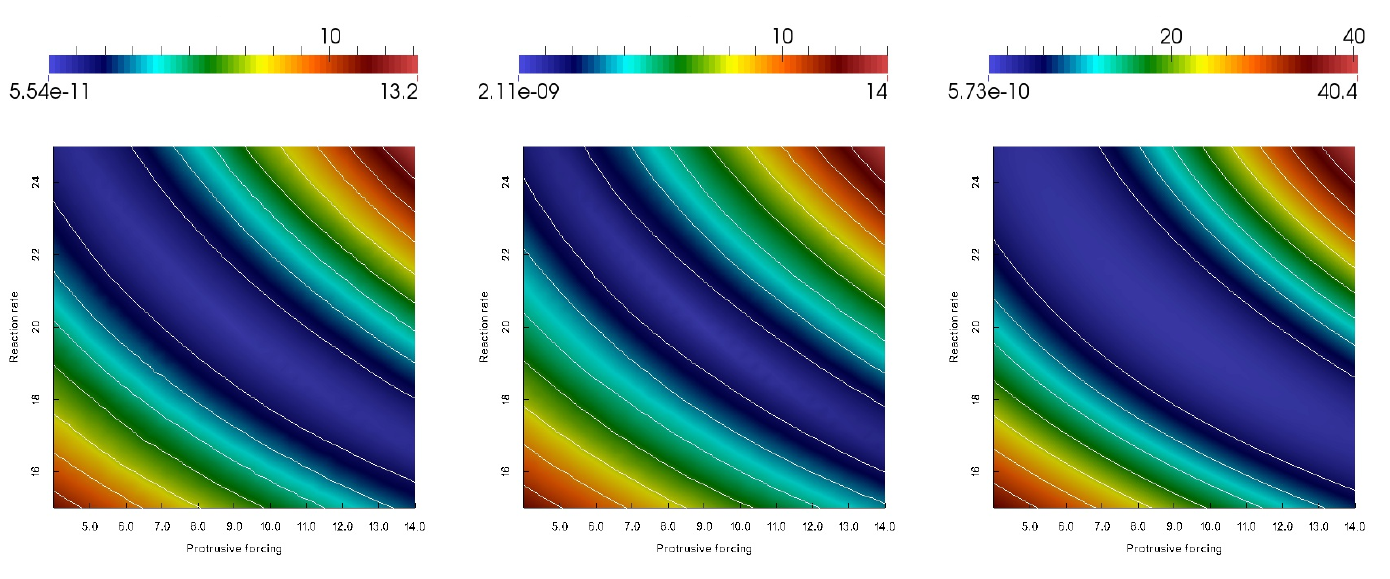}
\caption{(Colour online) The phase field objective functional (\ref{eqn:obj_pf}) as a function of the reaction rate $\gamma$ and the  protrusive forcing strength $(k_p)_2$ for phase field interfacial width parameter $\epsilon=0.001,0.01$ and $0.1$ reading from left to right. The observations were generated using $\gamma=20$ and $(k_p)_2=0.05$. The lines indicate level curves of the objective functional.  }\label{fig:objective_epsilon}
\end{figure}

Figure \ref{fig:objective_epsilon} shows the dependence of the phase field formulation of the objective functional on $\gamma$ and $(k_p)_2$ for different values of $\epsilon$ (the interfacial width parameter) and for all the weights set to one.  The sensitivity of the objective functional appears similar for the different values of $\epsilon$ and the values of the objective functional also appear to converge as we vary $\epsilon$ from $0.01$ to $0.001$ although the values do not  converge to the sharp interface formulation of the objective functional (Figure \ref{fig:objective} left hand column middle). This is unsurprising as the Haussdorff distance is used in the sharp interface setting and the $\Lp{2}$ distance is used in the phase field setting. Further investigations are warranted into the sharp interface limit of the phase field formulation of the objective functional.

\subsection{Numerical experiments for surfaces}\label{subsec:surfaces} 
We now consider the case where the surface  $\G\subset\Reals^3$. We generate the target data by approximating the solution to Problem \ref{pbm:forward} with a forcing term of the form (\ref{eqn:prop_forcing}) and reaction kinetics (\ref{eqn:schnak}). We selected  parameter values for the RDS and surface evolution law stated in Table \ref{tab:surface_parameter_values}. We approximated the solution to Problem \ref{pbm:forward} using linear finite elements on a mesh with 1026 degrees of freedom and selected a timestep of $10^{-2}$.
\begin{table}
\centering
\begin{tabular}{cccccccccc}
\toprule
$T$&$\sigma$&$k_b$&$D_1$&$D_2$&$(k_p)_1$&$(k_p)_2$&$\gamma$&$k_1$&$k_2$\\
\midrule
$5$&$5\times10^{-3}$&$0$&$1$&$100$&$-1\times10^{-2}$&$5\times10^{-2}$&$20$&$0.1$&$0.9$\\
\bottomrule
\end{tabular}\\
\caption[]{Parameter values used to generate the target data for numerical experiments on surfaces.}\label{tab:surface_parameter_values}
\end{table}

For the identification problem we attempted to recover the coupling term $(k_p)_2$ appearing in the forcing function $g$,  the parameter $\gamma$ that appears in the reaction kinetics (\ref{eqn:schnak}) and the surface tension of the cell membrane $\sigma$.  We assumed the target data was observed at $t=0,1,2,\dots,5$ and thus had $5$ observations of the data to fit to and one observation for the initial data. For the phase field  objective functional we computed the objective functionals on cubical domains such that the distance between the computed and target curves and the boundary of the domain was at least $2\epsilon$ and we used a triangulation of the cubes with 22065 degrees of freedom. 

Table \ref{tab:surface_identification} shows the results of two experiments, one for each choice of the objective functional (\ref{eqn:disc_obj_sharp}) or (\ref{eqn:obj_pf}). The starting values for the parameters, used as an initial guess for the algorithm, together with the computed values are reported together with the relative percentage errors.  The algorithm converged in 8 and 26 iterations for the sharp interface and phase field objective functionals respectively and took 2843 and 40217 seconds of CPU time for the sharp interface and phase field objective functionals respectively. The algorithm appears to perform best in the sharp interface case with the relative errors of two of  the estimated parameters in the phase field case exceeding the relative errors of the initial guesses.  In both cases the algorithm exhibits the smallest relative error in the identification of the reaction rate $\gamma$, with the estimates of the protrusive forcing strength $(k_p)_2$ slightly worse and the  with the relative error for the estimated value of the surface tension $\sigma$ significantly higher.

\begin{table}
\centering
\begin{tabular}{ccccc}
\toprule
&Parameter&True value&Starting value &Computed value\\
\midrule
\multirow{3}{0.25\textwidth}{Sharp interface objective functional: $\mmm J_{sharp}$}
&$(k_p)_2$&$5\times10^{-2}$&$3.75\times10^{-2}\ (25)$&$5.10\times10^{-2} \  (2.12)$\\
&$\gamma$&$20.00$&$25.00 \ (25)$&$19.77 \  (1.13)$\\
&$\sigma$&$5\times10^{-3}$&$5.75\times10^{-3}\ (15)$&$4.34\times10^{-3} \  (13.15)$\\
\\
\multirow{3}{0.25\textwidth}{Phase field objective functional: $\mmm J_{pf,\epsilon}$}
&$(k_p)_2$&$5\times10^{-2}$&$3.75\times10^{-2}\ (25)$&$3.70\times10^{-2} \  (25.90)$\\
&$\gamma$&$20.00$&$25.00 \ (25)$&$20.90 \  (4.51)$\\
&$\sigma$&$5\times10^{-3}$&$5.75\times10^{-3}\ (15)$&$1.82\times10^{-4} \  (96.37)$\\
\bottomrule
\end{tabular}\\
\caption[]{Results for an identification experiment for the forced motion of a surface, relative $\%$ error in parentheses.}\label{tab:surface_identification}
\end{table}             

\begin{figure}[ht!]
\centering
\includegraphics[trim = 0mm 0mm 0mm 0mm, clip=true, width=\linewidth]{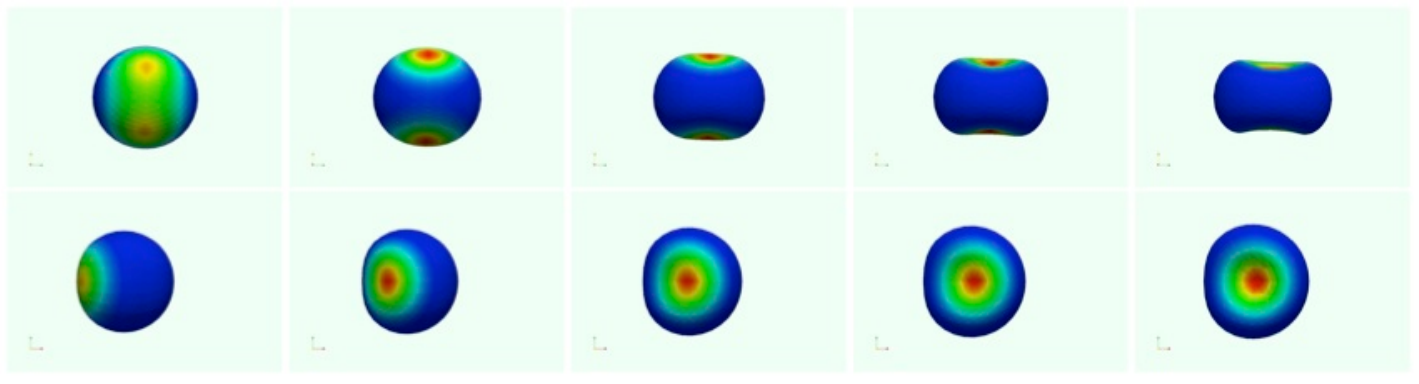}
\caption{(Colour online) Snapshots of the computed surface for the simulations described in \S \ref{subsec:surfaces} shaded by activator concentration $a_1$ at $t=1,2,3,4$ and $5$, from two different viewpoints computed using the starting parameter values given in Table \ref{tab:surface_identification}, for remaining parameter values see text.}\label{fig:surface_starting}  
\end{figure}

\begin{figure}[ht!]
\centering
\includegraphics[trim = 0mm 0mm 0mm 0mm, clip=true, width=\linewidth]{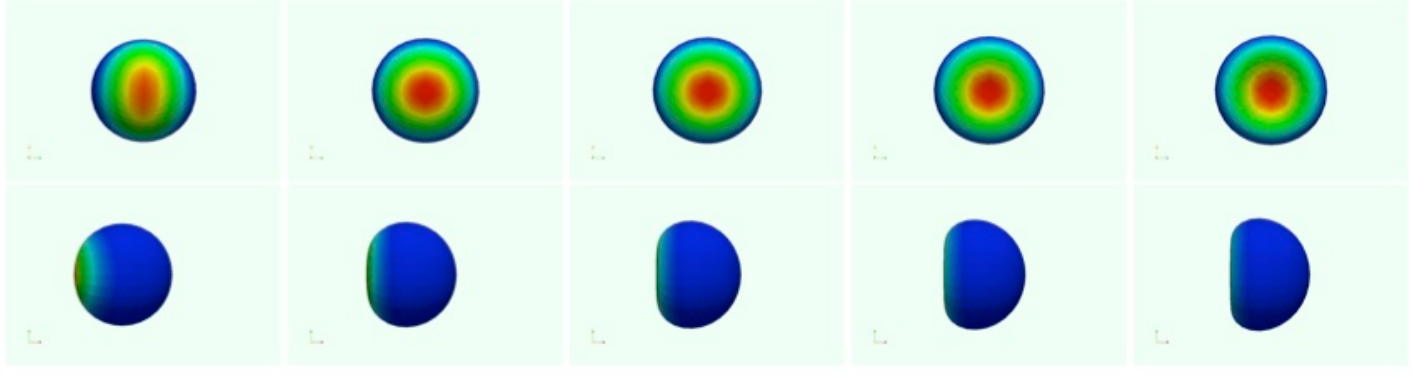}
\caption{(Colour online) Snapshots of the computed surface for the simulations described in \S \ref{subsec:surfaces}  shaded by activator concentration $a_1$ at $t=1,2,3,4$ and $5$, from the same viewpoints as in Figure \ref{fig:surface_starting}. The snapshots were computed using the estimated parameter values corresponding to the sharp interface objective functional given in Table \ref{tab:surface_identification}.}\label{fig:surface_optimal}  
\end{figure}

Figure \ref{fig:surface_starting} shows snapshots of the computed surface shaded by activator concentrations computed at the starting values used in the identification algorithm and Figure \ref{fig:surface_optimal} shows snapshots of the computed surface shaded by activator concentrations computed at the estimated parameter values computed by the identification algorithm using the sharp interface objective functional. We observe that qualitatively the solutions using the starting values which serve as an initial guess for the algorithm are far from the solutions using the estimated values, with the starting values generating a ``2-peak'' solution on the surface and markedly different motion. Despite this the proposed algorithm is sufficiently robust that  the results obtained with the estimated parameter values are qualitatively indistinguishable from those obtained with the true parameter values (not reported).

\begin{figure}[ht!]
\centering
\includegraphics[trim = 0mm 0mm 0mm 0mm, clip=true, width=\linewidth]{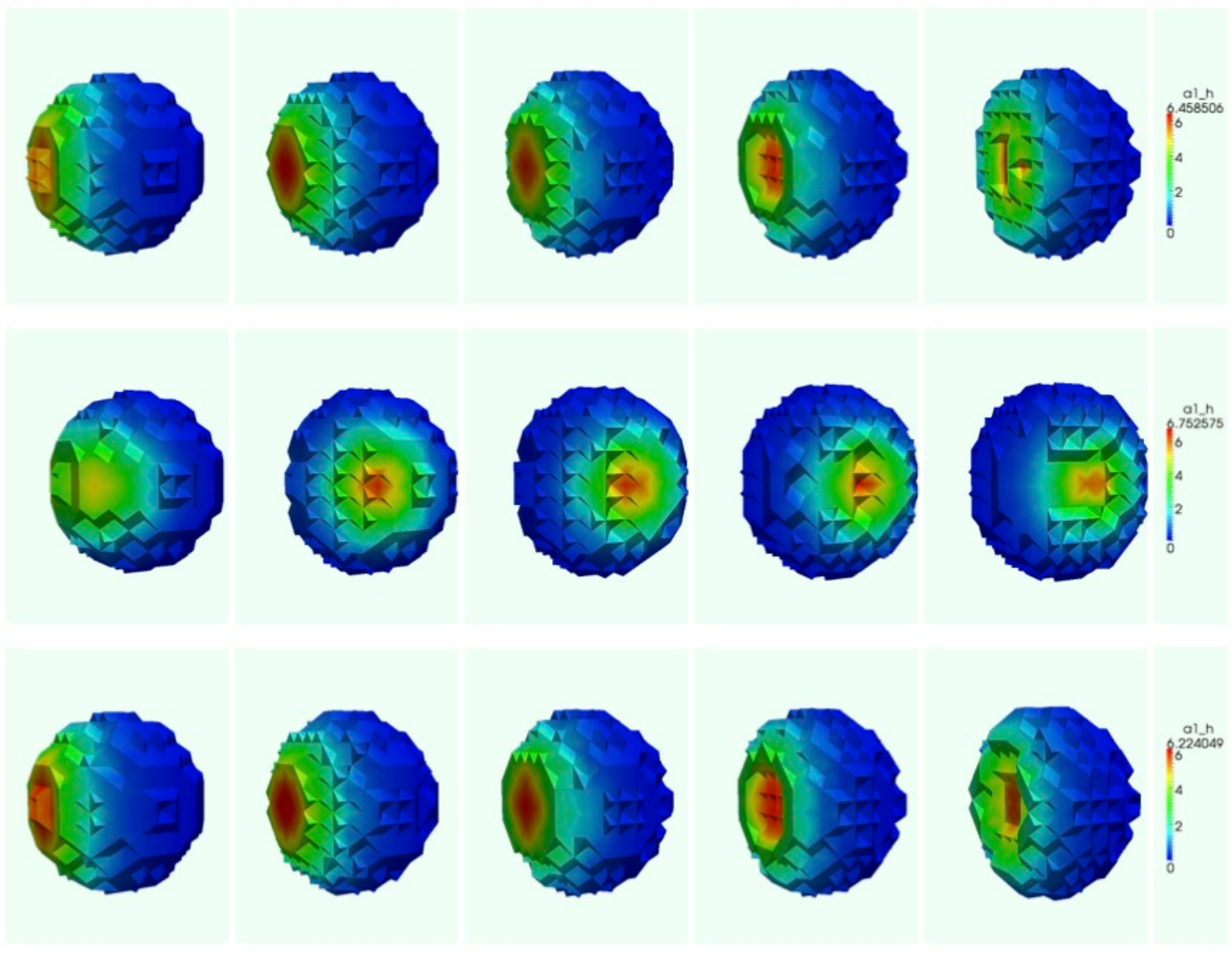}
\caption{(Colour online) Snapshots of the $-0.5$ level set of the phase field representation of the surfaces for the simulations described in \S \ref{subsec:surfaces}. The functions $\phi_\epsilon$ and $\hat{\phi}_\epsilon$ (c.f., (\ref{eqn:phi}) and (\ref{eqn:hat_phi}))  shaded by the phase field representation of the activator concentration $a_1$ i.e, the functions $(a_1)_\epsilon$ and $(\hat a_1)_\epsilon$ are shown at $t=1,2,3,4$ and $5$ (reading from left to right in each row). The top row corresponds to the phase field representation of the observations, the middle row to the simulated results using the initial guesses for the parameter values and the bottom row the simulated results with the estimated parameter values.}\label{fig:surface_phasefield}  
\end{figure}

Figure \ref{fig:surface_phasefield} shows snapshots of the $-0.5$ level set of the phase field representation of the surfaces shaded by the phase field formulation of the activator concentration on the surfaces. The top row of the figure corresponds to the observations, i.e., the data computed using the true parameter values, the middle row corresponds to the data computed using the initial guess for the parameter values and the bottom row the data computed using the estimated parameter values. We observe that while the initial guess for the parameter values generates a solution that is qualitatively dissimilar to the solution with the true parameter values, the solutions with the estimated parameters are qualitatively similar to the solutions with the true parameter values. This suggests that there may be multiple local minima of the objective functional in the phase field setting which may account for the large relative errors in the estimates of the parameters. Indeed one may expect  the objective functional to be relatively insensitive to positively correlated perturbations  in the surface tension (primarily a contractive force as the mean curvature of the surface is generally positive in this example) and the protrusive forcing strength $(k_p)_2$, which may account for the estimates for both of these parameters being lower than the true value.


\section{Numerical experiments with real data}\label{sec:real_data}
We now investigate the performance of the identification approach presented in this study when the observations consist of experimental rather than artificial data. We focus on modelling monopolar growth of fission yeast {\it Schizosacchromyces pombe}; a rod-shaped organism that proliferates asexually through growth at cell poles \citep{bendezu2012cdc42}. Following division, cells begin growth at a single tip (the tip that existed prior to division) before a minimal length is achieved enabling transition to bipolar growth. Bipolar growth continues until cell division, when growth machinery is relocated to the septum in the centre of the cell \citep{kelly2011spatial}. It has become apparent that, the conserved small Rho-like GTPase Cdc42, plays an essential role in regulating cell this polar cell growth. Accumulation of active Cdc42 at a growth tip, defines an area where vesicle delivery, exocytosis and cell wall remodelling occurs thereby delivering cell wall synthases to promote growth \citep{bendezu2012cdc42,Das13072012,drake2013model}. 

\subsection{Experimental observations of the monopolar growth of fission yeast cells}
\begin{figure}[ht!]
\centering
\includegraphics[trim = 0mm 0mm 0mm 0mm, clip=true, width=\linewidth]{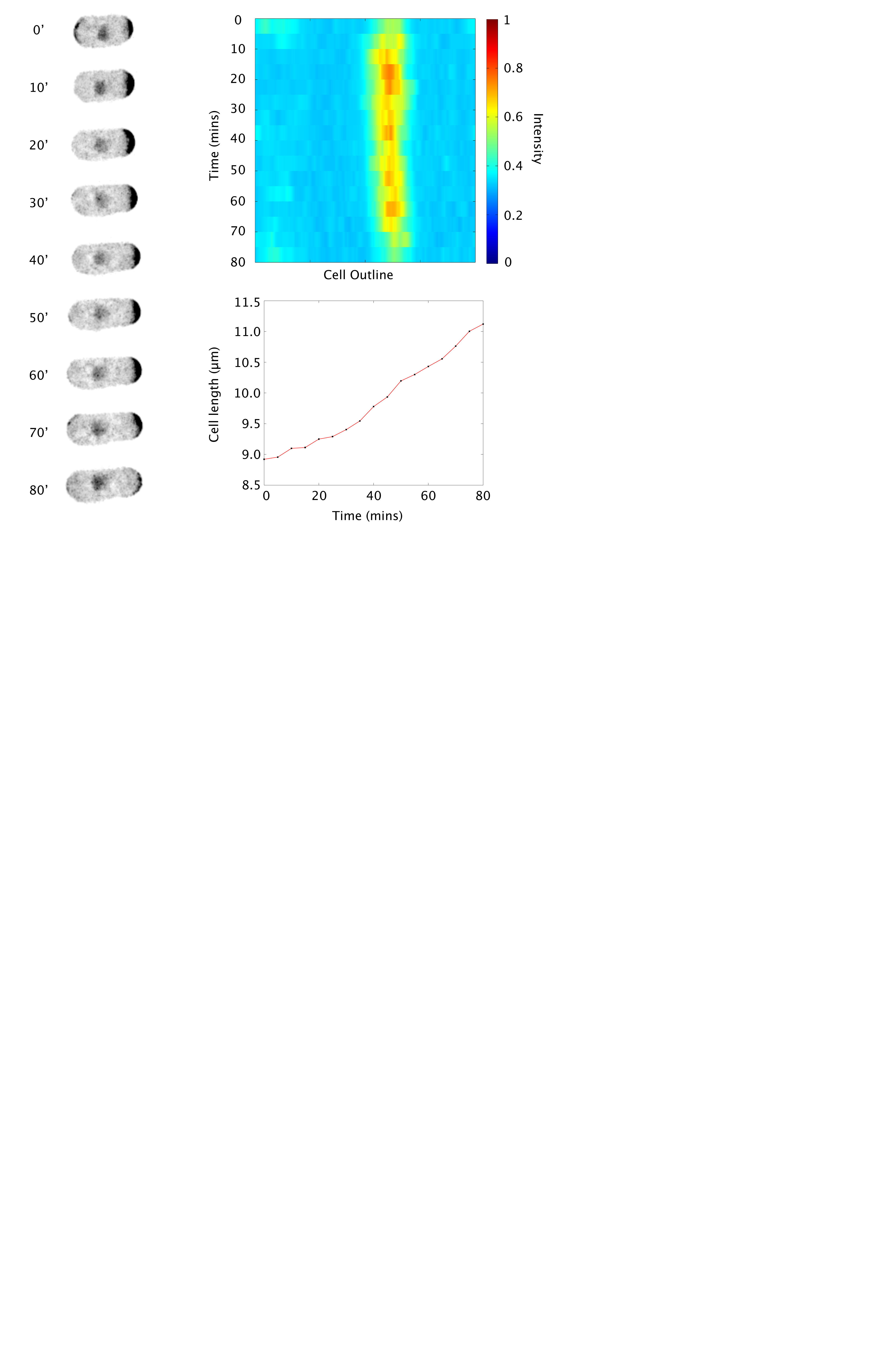}
\caption{(Colour online) Spatiotemporal profile of Cdc42 activity and cell growth in fission yeast. Time-lapse images of CRIB-GFP fluorescence indicating Cdc42 activity at a single tip of a representative cell undergoing monopolar cell growth (left). The plasma membrane temporal fluorescence profile (upper right) and tip-to-tip cell lengths (lower right) were calculated upon segmentation of the time-lapse image.}\label{fig:exp_data}  
\end{figure}

Recently, the importance of active Cdc42-GTP in the coordination of {\it S. pombe} growth sites during mitotic growth has become apparent. In vivo it is possible to determine the cellular localisation of active Cdc42 through the expression of a GFP-tagged Cdc42/Rac-interactive binding (CRIB) domain \citep{tatebe2008pom1}. We \citep{weston2013coordination} and others \citep{tatebe2008pom1,bendezu2012cdc42,Das13072012,kelly2011spatial} have described the localisation of CRIB-GFP in a range of {\it S. pombe} strains. In mitotically growing cells CRIB-GFP localises predominately to poles of cells \citep{weston2013coordination}. Using a previously described strain (JY1645) expressing CRIB-GFP \citep{weston2013coordination}, we performed time-lapse fluorescence microscopy (cells were grown on rich media agarose plugs) and were imaged every 5 minutes for a period of 80 mins (Figure \ref{fig:exp_data}, left panel). We quantified CRIB-GFP fluorescence at the plasma membrane and cell growth using QuimP2 (Figure \ref{fig:exp_data}, upper and lower right panels)  \citep{bosgraaf2009analysis}  from the point of cell division until initiation of new-end take off (bipolar growth). Further details of the experimental setup are given in Appendix \ref{Ap:bio_methods}.

 We selected five cells from the population of ten cells which exhibited single tip growth (some of the cells did not exhibit any marked growth, others grew from both tips and others were in close proximity to each other and thus influenced each others movement).  We (manually) picked the longest interval in which single tip growth was  evident for each of the five cells to be the  observations for the identification algorithm. Therefore the data we attempted to fit to, shown in Figure \ref{fig:yeast_experiments} (top row),  consisted of five cells observed over 11, 7, 5, 9 and 11 snapshots respectively, i.e., 55, 35, 25, 45 and 55 minutes respectively. For each of the cells, the initial snapshot at the start of the period of sustained monopolar growth was used to define the initial data for the model equations. 

\subsection{Modelling the monopolar growth of fission yeast cells}
 Experimental evidence suggests that growth is positively correlated with Cdc42 concentration and that there is a threshold value of Cdc42 concentration above which growth occurs \citep{Das13072012}.  In contrast to  \citet{drake2013model}, where a model with a Cdc42 signal distributed over a characteristic length-scale rather than generated by a reaction-diffusion process is employed, we propose a simple model for the Cdc42-dependent monopolar growth of fission yeast cells of the form considered in Problem  \ref{pbm:forward}. We propose that for the short timescale monopolar growth  that we seek to model, the evolution of the surface Cdc42 concentration satisfies a surface heat equation i.e., in place of (\ref{eqn:RDS}) we consider a one species heat equation with no reaction kinetics. Furthermore, as Cdc42 remains elevated and localised at the new tip we select a small membrane diffusion coefficient. This surface heat equation is coupled to the evolution law (\ref{eqn:evolution_law}). The volume of the cells is not conserved as they grow and thus the penalisation term $\lambda$ appearing in (\ref{eqn:evolution_law}) is set to zero.  Observations suggest that surface tension and bending rigidity only play minor roles for the short timescale monopolar growth considered in this study therefore we select small positive values for the surface tension and the bending rigidity. The parameter values used in the simulations are given in Table \ref{tab:yeast_parameter_values}.

 \begin{table}
\centering
\begin{tabular}{ccc}
\toprule
$\sigma$&$k_b$&$D$\\
\midrule
$1\times10^{-2}$&$1\times10^{-2}$&$1\times10^{-3}$\\
\bottomrule
\end{tabular}\\
\caption[]{Fixed parameter values used in the simulations of monopolar growth of fission yeast cells}\label{tab:yeast_parameter_values}
\end{table}

  To approximate the experimentally observed dependence of growth on Cdc42, i.e., growth above a given threshold value,  we consider a forcing function ($g$ in (\ref{eqn:evolution_law})) of the form
 \begin{equation}\label{eqn:yeast_forcing}
 g(\eta)=
 \begin{cases}
 0\quad &\mbox{if }\eta<k_1,\\
 k_2\left(\left(\frac{(\eta-k_1)}{k_{reg}}\right)\left(3-\frac{2\left(\eta-k_1\right)}{k_{reg}}\right)\right)^2\quad&\mbox{if } k_1<\eta<k_1+k_{reg},\\
 k_2\quad&\mbox{otherwise}, 
 \end{cases}
 \end{equation}
  where the parameters $k_1,k_2\geq0$ are the parameters we seek to identify. The parameter $k_{reg}\geq 0$ is a small positive parameter which governs the width of transition region between no growth and maximal growth. Threshold dependence of growth on Cdc42 is approximated as the parameter $k_{reg}$ approaches zero.  Although in theory we could attempt to identify the size of the transition region $k_{reg}$,  our numerical studies suggest that there are insufficient data points within this threshold region in each of the snapshots to make robust identification of this parameter possible (i.e., the spatial resolution of the data generated by the cell segmentation algorithm is too low). We therefore select $k_{reg}=5\times{10}^{-2}$ in all the subsequent simulations. The parameters we seek to identify may be interpreted as follows,
  $k_1$ is the lower threshold below which there is no growth, $k_1+k_{reg}$ is the upper threshold above which growth is maximal and $k_2$ the maximal strength of  the Cdc42 dependent forcing. 
\begin{Rem}[Modelling simplifications]
The mechanism that governs localisation of Cdc42 at the new tip is complex and involves many other species such as Gef1 \citep{Das13072012}, Ras1-GDP, Ras1-GTP, Gap1 \citep{weston2013coordination} and Scd1 \citep{onken2006compartmentalized}. Modelling these complex multi-species interactions are beyond the scope of this study hence we consider the simplified  model outlined above.
 
 The inclusion of a small positive surface diffusion coefficient for Cdc42 and small positive bending rigidity and surface tension  serve to regularise the surface PDE - surface evolution law system making them suitable for approximation with the surface finite element method utilised in this study. We have verified that doubling either  the diffusion coefficient, bending rigidity or surface tension yields relatively unchanged parameter estimates.
 \end{Rem}
 \begin{Rem}[Comparison of the identified parameters with those estimated in the literature]
 Due to the structure of the proposed evolution law modelling monopolar growth, in which the forcing due to Cdc42 is the dominant term, the parameter $k_2$  may be  compared with experimental data on the growth rates of tips with elevated Cdc42 (specifically, the old tip that existed prior to division)  as reported in \citep[Fig. 1F]{Das13072012}. We therefore report on our estimates of this parameter in units of $\mu m/min$ to allow comparison with the experimental measurements of \citet{Das13072012}.
 
 On the other hand, the threshold level $k_1$ may be thought of as a cell specific quantity and it is not possible to directly relate this to a concentration level of Cdc42. This is due to the fluorescence source, CRIB-GFP binding to Cdc42-GTP (the active conformation). The increased fluorescence observed, for example at a growing tip, is due to CRIB-GFP accumulation upon its binding to active Cdc42. Whole cell fluorescence intensity would reflect the total cellular concentration of CRIB-GFP but not Cdc42. The concentration of Cdc42 at the plasma membrane cannot be directly inferred by fluorescence intensity measurements as fluorescence is dependent on both the proportion of Cdc42 that has become activated and the concentration of CRIB-GFP available for binding. Variability in expression levels from cell-to-cell will cause variation in overall Cdc42 and CRIB-GFP concentration, therefore contributing to cell-to-cell variability in the levels of fluorescence intensity that are observed.

 \end{Rem}

  \subsection{Identification experiments for the monopolar growth of fission yeast cells}
\begin{figure}[ht!]
\centering
\includegraphics[trim = 0mm 0mm 0mm 0mm, clip=true, width=\linewidth]{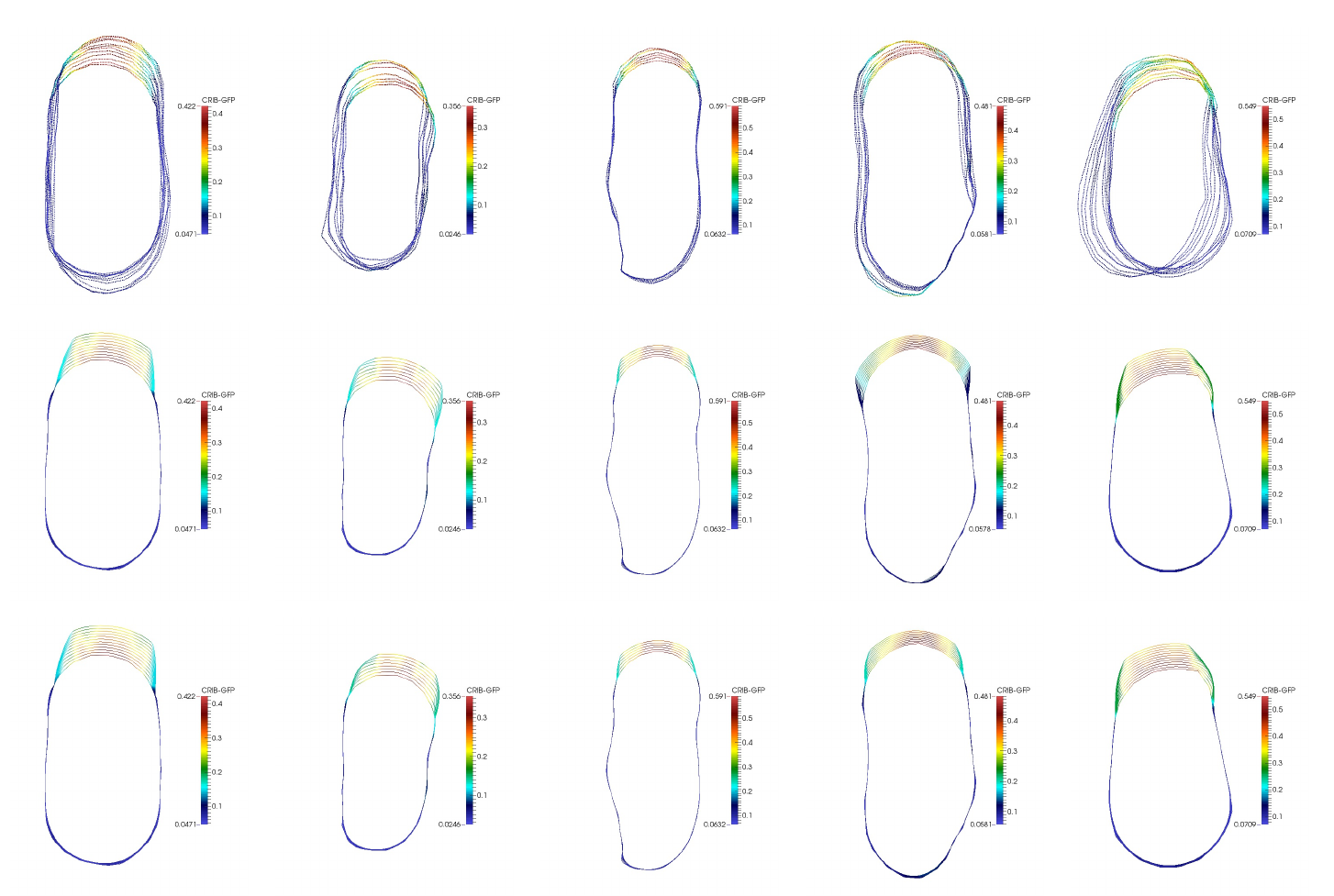}
\caption{(Colour online) Monopolar growth of five yeast cells shaded by CRIB-GFP intensity. Top row: experimental observations, middle row: solution with identified parameters using the sharp interface formulation of the objective functional with equal weights and bottom row: solution with identified parameters using the phase field formulation of the objective functional with equal weights. For further details see text.
We observe that in all cases the identified parameters are such that the tip of the cell lacking Cdc42 does not deform significantly while the tip with elevated Cdc42 grows. The only cell for which the phase field and sharp interface formulations generate qualitatively different results is the cell in column four for which the lower threshold value identified in the sharp interface case, c.f., Figure \ref{fig:yeast_forcing_equal}, results in a larger proportion of the cell growing.}\label{fig:yeast_experiments}  
\end{figure}

\begin{figure}[hp!]
\centering
\subfigure[][Equal weights.]{
\includegraphics[trim = 0mm 0mm 0mm 0mm, clip=true, width=0.7\linewidth]{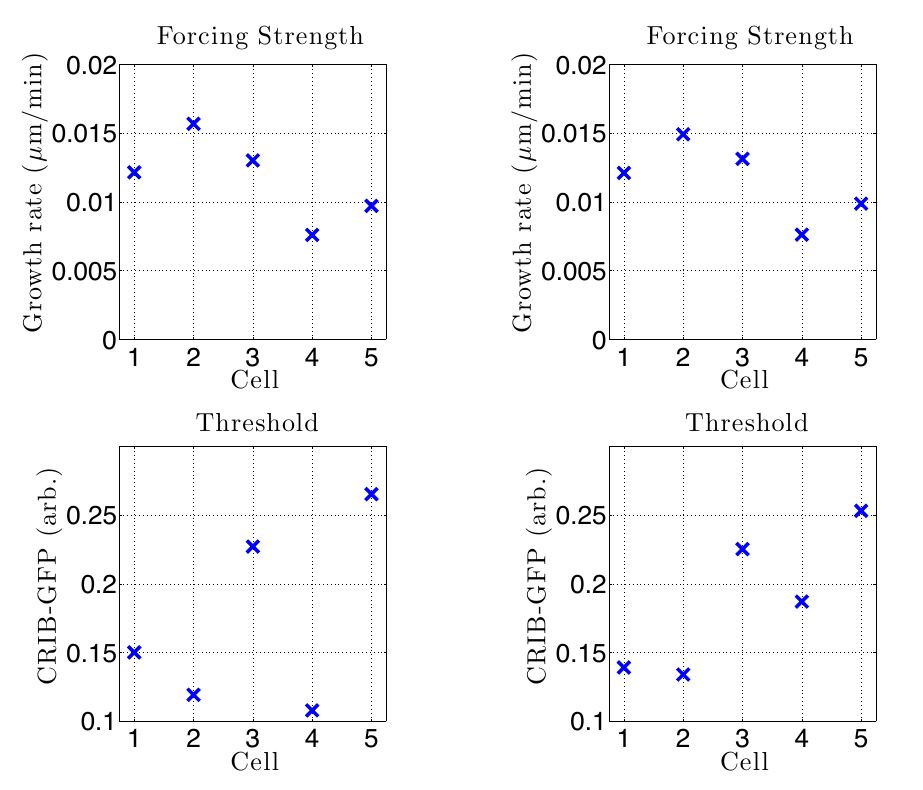}
\label{fig:yeast_forcing_equal}  
}
\subfigure[][Weights chosen such that the error due to position and concentration are comparable at initial parameter values.]{
\includegraphics[trim = 0mm 0mm 0mm 0mm, clip=true, width=0.7\linewidth]{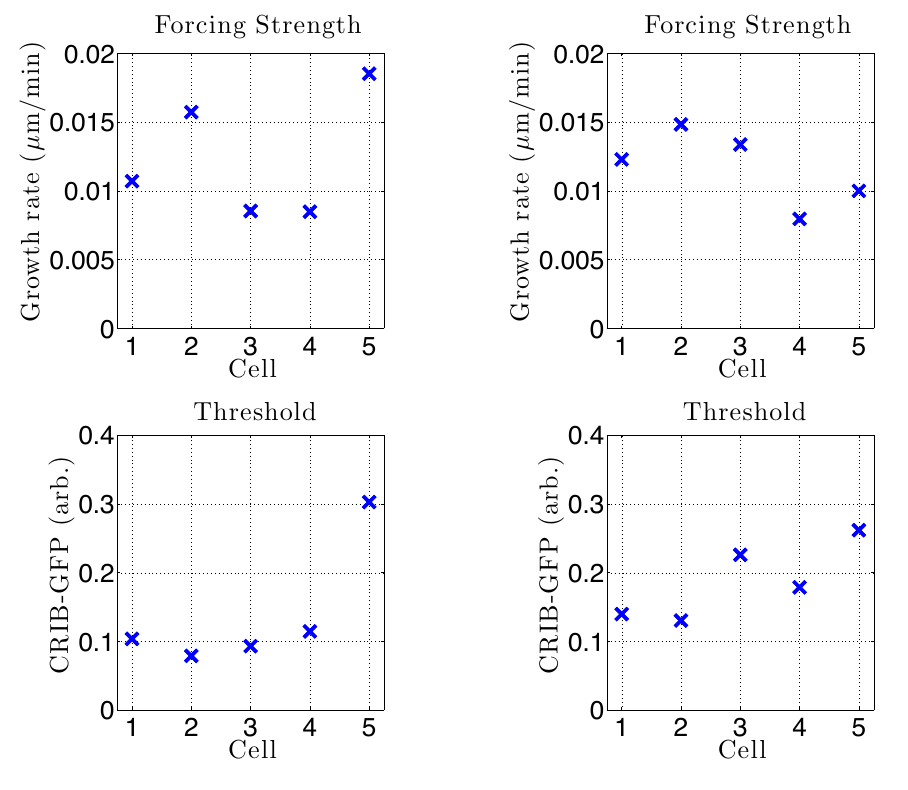}
\label{fig:yeast_forcing_weighted}  
}
\caption{(Colour online) The identified parameters for the numerical experiments described in \S \ref{sec:real_data}. Within each figure the top row shows the estimates of the maximal forcing strength $k_2$ and the bottom row the threshold intensity $k_1$ (c.f., (\ref{eqn:yeast_forcing})), the left column corresponds to the sharp interface objective functional and the right column the phase field objective functional. The ordering of the cells within each plot is as in Figure \ref{fig:yeast_experiments}.}\label{fig:yeast_forcing}  
\end{figure}

  For all the simulations in this section the value of the interfacial width parameter $\epsilon$ related to the phase field version of the algorithm, c.f., (\ref{eqn:obj_pf}), was taken to be  $0.63\mu m$. The remaining parameters used in the LM algorithm were the same as those employed in \S \ref{sec:examples}.  For the first identification experiment we set each of the weights to be one (i.e., $w_i=1,i=1,\dots,2n_s$ in (\ref{eqn:chi_sharp}) and (\ref{eqn:chi_pf})) and selected initial parameter estimates of $0.1$ (arbitrary units) for the threshold value $k_1$ and $0.25 \mu m/min$ for the maximal forcing strength $k_2$. Figure \ref{fig:yeast_experiments} shows experimental observations of the monopolar growth of the four cells considered over a three minute interval with the cells shaded by fluorescence intensity of Cdc42 (top row) and the solution with the identified parameters using the sharp interface and phase field objective functionals (middle and bottom rows respectively).   The parameters identified by the algorithm with the two different objective functionals are shown in Figure \ref{fig:yeast_forcing_equal}  for the sharp interface (left) and phase field (right) objective functionals. The algorithm selects a  threshold  $(k_1 \text{ in } (\ref{eqn:yeast_forcing}))$  of between 0.1 and 0.3 for the CRIB-GFP intensity below which there  is no growth. The maximal forcing strength $(k_2 \text{ in } (\ref{eqn:yeast_forcing}))$ identified is between 
  0.005 and 0.02 and this is in agreement with the values measured experimentally in \citet[Fig. 1F]{Das13072012} . The identified parameters are very similar for both the different objective functionals with virtually identical forcing strength estimates and similar threshold estimates. Only the threshold ($k_1$) parameter identified for cell number 4 appears markedly different for the two objective functionals with the sharp interface functional selecting a significantly lower threshold value. The simulated results with the identified parameters Figure \ref{fig:yeast_experiments} bottom two rows, shows that the algorithm identifies parameters that generate monopolar growth qualitatively similar to that observed in experiments with little to no growth at the tip where the fluorescence intensity of Cdc42 is low. There are only minor qualitative differences between the solution with the phase field and sharp interface objective functionals, for example for the cell shown in the fourth column of Figure \ref{fig:yeast_experiments}, we see that a larger region of the cell surface forms the growing tip and this is due to the significantly lower threshold $k_1$ identified by the sharp interface algorithm in this case, see Figure \ref{fig:yeast_forcing_equal}.

\begin{figure}[hp!]
\centering
\subfigure[][Equal weights.]{
\includegraphics[trim = 0mm 0mm 0mm 0mm, clip=true, width=\linewidth]{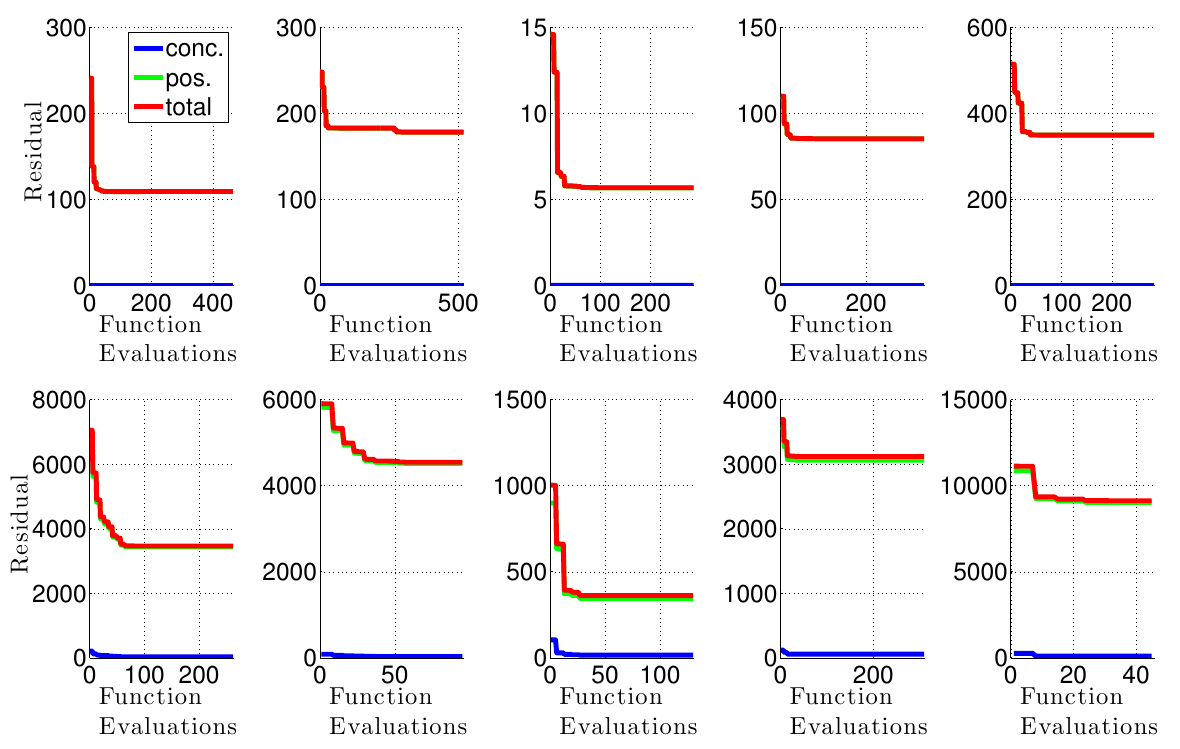}
\label{fig:yeast_residuals_equal}  
}
\subfigure[][Weights chosen such that the error due to position and concentration are comparable at initial parameter values.]{
\includegraphics[trim = 0mm 0mm 0mm 0mm, clip=true, width=\linewidth]{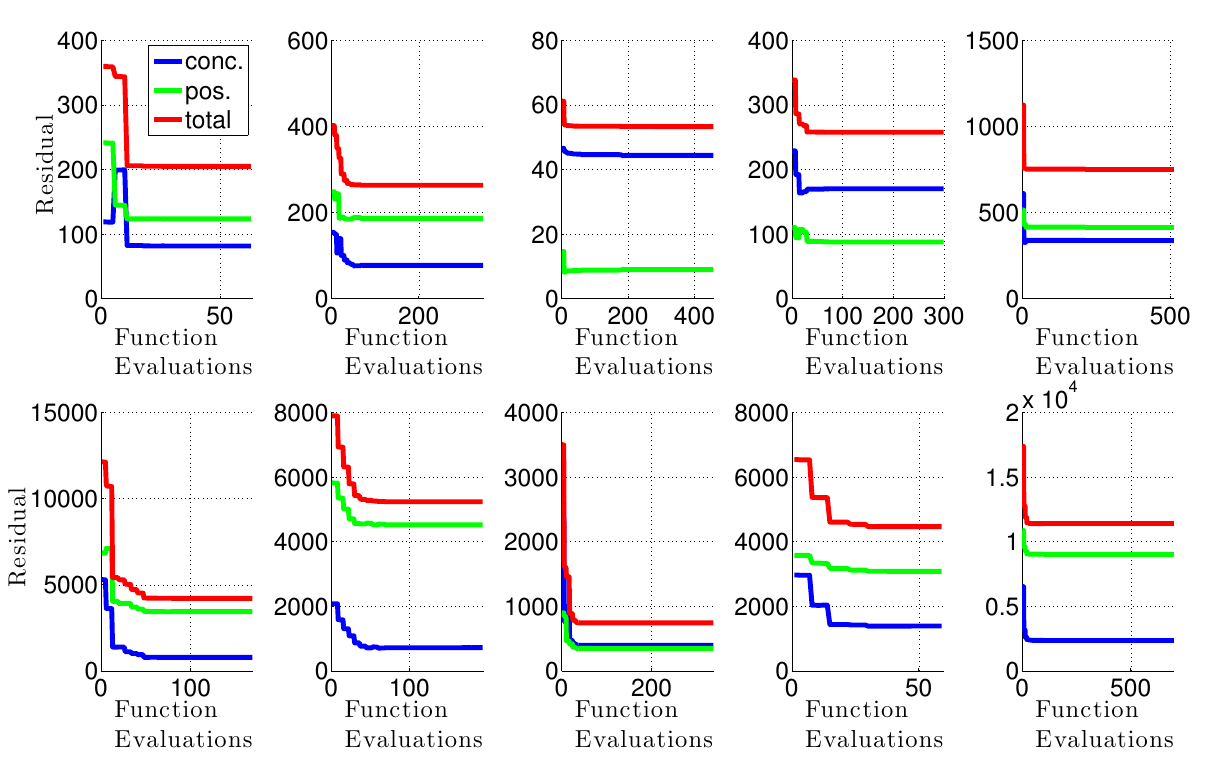}
\label{fig:yeast_residuals_weighted}  
}
\caption{(Colour online) The residual $J$ and the contribution of the the error due to position and the error due to concentration for the numerical experiments described in \S \ref{sec:real_data}. Within each figure the top row corresponds to the sharp interface objective functional and the bottom row the phase field objective functional. The ordering of the columns within each subfigure is as in Figure \ref{fig:yeast_experiments}. (For further details see text).}\label{fig:yeast_residuals}  
\end{figure}

Figure \ref{fig:yeast_residuals_equal} shows the evolution of the residual versus the number of function evaluations performed by the algorithm  for the sharp interface (top row) and phase field (bottom row) objective functionals. We observe the behaviour typical of the LM algorithm of initially rapid reduction of the residual followed by periods of stagnation.
In four of the five cases considered the phase field version of the algorithm converges with less than half the function evaluations needed by the sharp interface version of the algorithm, however, in the fourth column we see that the sharp interface version converges after a similar number of function evaluations as the phase field version. We also report on the component of the residual due to the error in position  ($\sum_{i=1}^{n_s}\chi_i^2$) and the component of the residual due to the error in concentration ($\sum_{i=n_s+1}^{2n_s}\chi_i^2$), where we have used the notation of \S \ref{sec:optimisation}. We observe that the error due to position dominates the error due to concentration for both formulations. 

\begin{figure}[ht!]
\centering
\subfigure[][Equal weights, sharp interface objective functional.]{
\includegraphics[trim = 0mm 0mm 0mm 0mm, clip=true, width=0.45\linewidth]{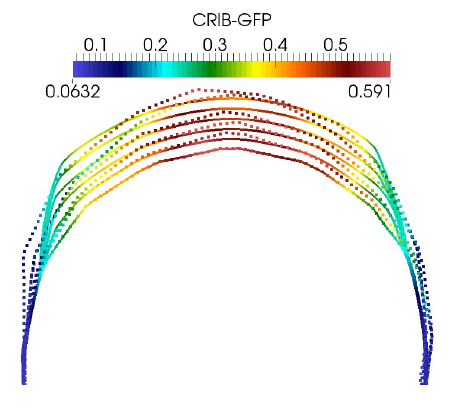}
\label{fig:yeast_tip_equal}  
}
\subfigure[][Weights chosen such that the error due to position and concentration are comparable at initial parameter values, sharp interface objective functional.]{
\includegraphics[trim = 0mm 0mm 0mm 0mm, clip=true, width=0.48\linewidth]{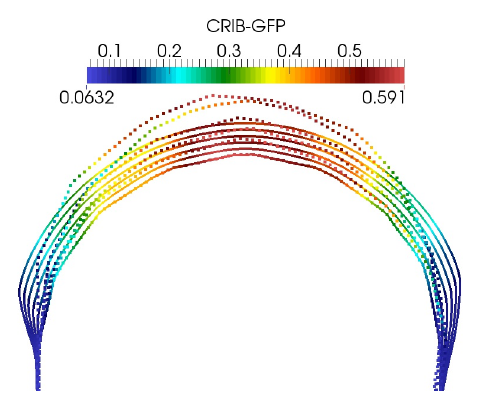}
\label{fig:yeast_tip_weighted}  
}
\subfigure[][Equal weights, phase field  objective functional.]{
\includegraphics[trim = 0mm 0mm 0mm 0mm, clip=true, width=0.45\linewidth]{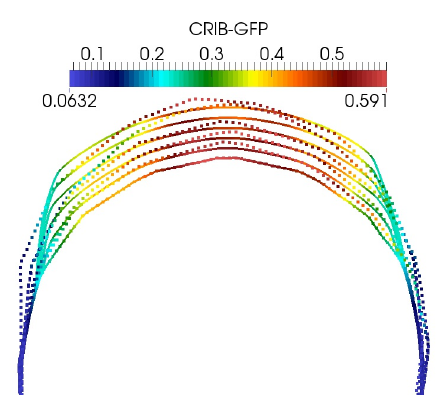}
\label{fig:yeast_tip_equal_pf}  
}
\subfigure[][Weights chosen such that the error due to position and concentration are comparable at initial parameter values, phase field  objective functional.]{
\includegraphics[trim = 0mm 0mm 0mm 0mm, clip=true, width=0.45\linewidth]{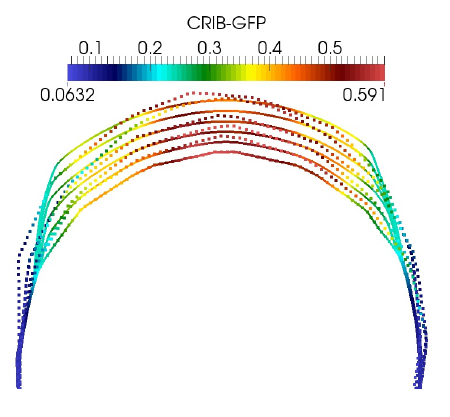}
\label{fig:yeast_tip_weighted_pf}  
}
\caption{(Colour online) Experimental data for the growing tip and snapshots of the solution with the identified parameters for the cell shown in the third  column of Figure \ref{fig:yeast_experiments}.}\label{fig:yeast_tips}  
\end{figure}

We conclude with another experiment where we selected weights such that at the initial guess for the parameter values the contribution of the error in concentration and position to the objective functional was roughly comparable. Specifically, we choose the weights such that the average over the five cells of the error due to position was the same as the error due to concentration at the initial guess for the algorithm which corresponded to  $w_i=1, w_{i+n_s}=5184, i=1,\dotsc,n_s$ for the sharp interface objective functional  and $w_i=1, w_{i+n_s}=25, i=1,\dotsc,n_s$ for the phase field objective functional. 
Figure \ref{fig:yeast_forcing_weighted} shows the  identified parameters for the sharp interface (right) and phase field (left) objective functionals with the new weights. In the sharp interface case (left column of Figure \ref{fig:yeast_forcing_weighted}) the forcing strength identified for cells 1,2 and 4 is similar to the equal weight case (Figure \ref{fig:yeast_forcing_equal}), for cells 3 and 5 different forcing strengths are identified by the algorithm, lower and higher respectively. For the phase field case the identified forcing is  very close to the parameters identified in the equal weights case for all five cells. The sharp interface algorithm identifies a lower threshold for cells 1 and 3 and a higher threshold for cell 5 than the equal weights case while in the phase field case the thresholds identified are similar to the  equal weights case for all five cells. Figure \ref{fig:yeast_residuals_weighted} shows the evolution of the residual versus  the number of function evaluations for the distinct weight case. There appears to be no clear relationship between the number of function evaluations necessary and the choice of either distinct or equal weights with some cells requiring more function evaluations and other cells less than in the equal weights case for both the sharp interface and phase field formulations of the objective functional.
  In the interest of space and as in the phase field case the results remain largely unchanged to those with equal weights, we do not report on all the simulations with the distinct weights, instead we illustrate the differences between the distinct weights and equal weights case by considering a single cell which the differences are evident. Figure \ref{fig:yeast_tips} shows magnifications of the data at the growing tip at a series of snapshots for one of the cells in Figure \ref{fig:yeast_experiments} together with the solution with the identified parameters at the corresponding times. There are differences in the solution with parameters identified by the sharp interface functional for the two different choices of weights. For the sharp interface objective functional, the equal weight case appears to fit, more accurately, the position of the growing tip while, the distinct weights case generates intensity profiles that reproduce the decline in CRIB-GFP intensity away from the growing tip. The phase field objective functional yields very similar solutions for both the equal and distinct weights, which is as expected since the parameter estimates shown in Figure \ref{fig:yeast_forcing} remain largely unchanged. This is in accordance with the results of \S \ref{subsec:sensitivity} where it was observed that the sharp interface formulation of the objective functional appears more sensitive to changes in weights than the phase field formulation of the objective functional. It also suggests that practitioners may tune the weights if the sharp interface formulation is used to prioritise fits to either the observed positions or the observed concentrations depending on the application in hand.

\section{Conclusion}\label{sec:conclusion}
The mathematical modelling of cell motility is a burgeoning research field \citep{hecht2011activated,herant2010cytopede,mogilner2009mathematics,ziebert2011model}. Typically it is difficult to estimate the forces associated with cell motility and it is therefore challenging to estimate parameters within these, often phenomenological, models. The advent of novel fluorescence microscopy, imaging and molecular biology techniques has meant   high resolution two and three dimensional data of migrating cells is  now available \citep{bosgraaf2009analysis} and this is one resource which may be used to estimate these parameters.       

In \citet{venk11chemotaxis} a general framework for modelling cell motility was proposed coupling a geometric evolution law for the motion of the cell membrane to partial differential equations posed on the membrane which incorporated some of the existing models for cell motility \citep{neilson2010use,neilson2011chemotaxis,neilson2011modelling}. In this study we proposed an algorithm for the identification of parameters in models  that fit into the general framework proposed in  \citet{venk11chemotaxis}, that makes use of the robust and efficient numerical method proposed in the same study for the approximation of the model equations. To our knowledge this is one of the first algorithms that allows parameters in models for cell motility to be estimated from experimental imaging data. We formulated the identification problem as a minimisation problem where the objective functional to be minimised may be computed either within a sharp interface or a phase field setting and presented an optimisation method based on the widely used Levenberg-Marquardt algorithm for the solution of the minimisation problem. We discussed in detail the implementation of the algorithm in this context and presented a number of numerical experiments where the algorithm was used to identify parameters with artificial data generated by simulating the model equations. The numerical tests indicate the algorithm is robust to moderate levels of noise and is capable of identifying parameters in three-dimensional models for cell motility and in the presence of large deformations. On the evidence of the numerical tests carried out, it appears that the sensitivities of the two different formulations of the objective functional exhibit significant differences. The algorithm equipped with the sharp interface formulation of the objective functional gives better estimates of certain parameters and is less computationally intensive, while the algorithm equipped with the phase field formulation of the objective functional gives estimates of certain parameters which are more robust to noise and  may be more amenable to analysis, which is an area we intend to explore in future work. We finished with an example where we illustrated the performance of the algorithm when the data consists of real experimental observations of cells migrating in vitro. The algorithm was used to identify parameters in a simple model for the monopolar growth of fission yeast cells. The algorithm identifies threshold values and maximal growth rates under a model for growth where growth occurs above  threshold value of Cdc42 concentration, the maximal growth rates identified  agree with the related maximal growth rate of tips  measured experimentally in \citet{Das13072012}. Furthermore, by tuning the weights in the sharp interface version of the objective functional it appears one may prioritise the fit to the observed concentrations or to the observed positions. Thus there is some flexibility for practitioners dependent on the application in hand.

We hope the proposed algorithm  will be used as modelling tool. For example the algorithm may be used to identify which parameters are most likely to have been changed  in a mutant  that exhibits qualitatively different motility. We also intend to consider other optimisation methods such as adjoint based methods arising from optimal control theory and Bayesian methods and to compare and contrast such methods with the method proposed in this study.

\section*{Acknowledgments}
This work was supported by a University of Warwick Impact Fund grant (C.V, C.M.E. and B.S.), the Engineering and Physical Sciences Research Council grants EP/G010404/1 (C.V and C.M.E.) and EP/J016780/1 (C.V.), a University of Warwick Impact Fund grant (C.W. and G.L.), the Biotechnology and Biosciences Research Council (W.C. and G.L.) grant number (BB/G01227X/1) and the Birmingham Science City Research Alliance (G.L.)


\appendix
\section{Biological methods}\label{Ap:bio_methods}
\subsection{Yeast strain and growth media}
The generation of the yeast strain (JY1645) used in this study is described elsewhere \citep{weston2013coordination}. General yeast procedures were performed as described previously \cite{davey199513,ladds1996sxa2} using rich amino acid (AA) media \citep{mos2013role}. 

\subsection{Microscopy}
Fluorescent time course experiments were performed on 2\% rich AA media agarose plugs at $29^\circ C$ with images acquired every 5 mins. CRIB-GFP was visualised using a Personal DeltaVision (Applied Precision, Issaquah, WA) comprising, an Olympus UPlanSApo 100x, N.A. 1.4, oil immersion objective and a Photometric CoolSNAP HQ camera (Roper Scientific). Captured images were processed by iterative constrained deconvolution using SoftWoRx (Applied Precession) and analyzed using ImageJ. Images used in analysis and representations are maximal projections of Z-stacks.

\subsection{Data analysis}
Image analysis was performed using the open source program ImageJ (http://rsb.info.nih.gov/ij/). Cell segmentation was achieved using the BOA plug-in component of the Quantitative Imaging of Membrane Proteins (QuimP) package (http://go.warwick.ac.uk/bretschneider/quimp) \citep{bosgraaf2009analysis,dormann2002simultaneous}. We have previously described the use of QuimP to quantify plasma membrane fluorescence intensities in fission yeast \citep{bond2013quantitative}.

 \end{document}